\documentclass{amsart}
\usepackage{amsmath,amssymb,amsthm,texdraw,citesort}

\numberwithin{equation}{section}

\theoremstyle{plain}
\newtheorem{theorem}{Theorem}[section]
\newtheorem{prop}[theorem]{Proposition}
\newtheorem{lemma}[theorem]{Lemma}
\newtheorem{cor}[theorem]{Corollary}
\theoremstyle{definition}
\newtheorem{definition}[theorem]{Definition}
\newtheorem{example}[theorem]{Example}
\theoremstyle{remark}
\newtheorem{remark}[theorem]{Remark}

\newcommand{\nc}{\newcommand}
\nc{\ol}{\overline} \nc{\eit}{\tilde{e}_i} \nc{\fit}{\tilde{f}_i}
\nc{\ali}{\alpha_i} \nc{\op}{\oplus} \nc{\ot}{\otimes}
\nc{\N}{\mathbf{N}} \nc{\Z}{\mathbf{Z}} \nc{\Q}{\mathbf{Q}}
\nc{\g}{\mathfrak{g}} \nc{\uq}{U_q} \nc{\ftil}{\tilde{f}}
\nc{\etil}{\tilde{e}} \nc{\an}{A_n} \nc{\anone}{A_n^{(1)}}
\nc{\B}{\mathcal{B}} \nc{\sln}{\mathfrak{sl}_{n+1}}
\nc{\slnh}{\hat{\mathfrak{sl}}_{n+1}}
\nc{\defi}[1]{\emph{\textbf{#1}}} \nc{\veps}{\varepsilon}
\nc{\vphi}{\varphi} \nc{\vepst}{\tilde{\varepsilon}}
\nc{\vphit}{\tilde{\varphi}} \nc{\La}{\Lambda} \nc{\la}{\lambda}
\nc{\T}{\mathcal{T}} \nc{\spaceosi}{\Phi} \nc{\M}{\mathcal{M}}
\nc{\Me}{\mathcal{M^{E}}} \nc{\Minfty}{\mathcal{M^{\infty}}}
\nc{\Mec}{\mathcal{M}_{c}^{\mathcal E}} \nc{\hwc}{\mathcal{B}}
\nc{\ra}{\rightarrow}

\nc{\Binf}{\mathcal{B}(\infty)} \nc{\Tinf}{\mathcal{T}(\infty)}

\DeclareMathOperator{\wt}{wt}
\DeclareMathOperator{\wtt}{\tilde{wt}}

\begin{document}

\title
[Monomial descriptions of crystals $\B(\infty)$ and $\B(\la)$ for
$A_n$] {Nakajima Monomials and crystals\\ for special Linear Lie
algebras}

\author{Hyeonmi Lee}
\address{Research Institute for Mathematical Sciences,
         Kyoto University, Kyoto, 606-8502 Japan}
\email{hyeonmi@kurims.kyoto-u.ac.jp, hmlee@kias.re.kr}
\thanks{This work was supported in part by KOSEF Grant
        R01-2003-000-10012-0}

\begin{abstract}
We present explicit descriptions of the crystals $\B(\infty)$ and
$\B(\la)$ over special linear Lie algebras in the language of
\emph{extended Nakajima monomials}. There is a natural
correspondence between the monomial description and Young tableau
realization, which is another realization of crystals $\B(\infty)$
and $\B(\la)$.
\end{abstract}

\maketitle

%%%%%%%%%%%%%%%%%%%%%%%%%%%%%%%%%%%%%%%%%%%%%%%%%%%%%%%%%%%%%%%%%%%%%%%%%%
\section{Introduction}%
\label{Introduction}

The theory of \emph{Nakajima monomials} is a combinatorial scheme
for realizing crystal bases of quantum groups. Nakajima introduced
a certain set of monomials realizing the irreducible highest
weight crystals in~\cite{MR1865400}. Kashiwara and Nakajima
independently defined a crystal structure on the set of Nakajima
monomials and also gave a realization of irreducible highest
weight crystal $\B(\la)$ in terms of Nakajima monomials, as the
connected component of the monomial set containing a maximal
vector of dominant integral weight
$\la$~\cite{MR1988989,MR1988990}. This has lead to the belief that
it should be possible to give a similar realization for
$\B(\infty)$, which is the crystal base of the negative part
$U_q^-(\mathfrak{g})$ of a quantum group over symmetrizable
Kac-Moody algebra $\g$, also.

Much effort has been
made~\cite{MR1475048,MR1286129,ks,saito,MR1614241,savage2,savage3,lee}
to give realization of $\B(\infty)$ over various Kac-Moody algebras.
In addition to these works, in our recent work~\cite{hm,HL}, we gave
new realization of $\B(\infty)$ for the finite simple Lie algebras,
in terms of Young tableaux.

Starting from the realization theorem of Kashiwara and
Nakajima~\cite{MR1988989,MR1988990}, we can argue that it is not
possible to find the crystal $\B(\infty)$ within the set of Nakajima
monomials with their given crystal structure. Hence, in our
work~\cite{lee}, we constructed the set of extended Nakajima
monomials and developed a crystal structure on it, and also gave
explicit descriptions of $\B(\infty)$ for $A_n^{(1)}$ case, in the
language of extended Nakajima monomials. Actually, the set of
Nakajima monomials can be embedded as a subcrystal in this set of
extended Nakajima monomials. Thus, the monomial theory developed for
irreducible highest weight crystal can easily be transferred to that
on the extended monomial set.

As the first contribution of our present paper, we introduce
explicit descriptions of the crystal $\B(\infty)$, in terms of
\emph{extended Nakajima monomials}. We restrict ourselves to special
Linear Lie algebras. The extended Nakajima monomial description is
obtained by relating it to the Young tableau realization~\cite{HL}.

The second contribution of this paper is to give an explicit
description of the irreducible highest weight crystal $\B(\la)$ for
any dominant integral weight $\la$, in monomial language, for $\an$
case. Another monomial description of $\B(\la)$ for this type may be
found in~\cite{MR2043368}, but unlike this work, there is an
immediate correspondence between our description and the Young
tableau realization of Kashiwara and Nakashima~\cite{MR1273277}.

Our paper is organized as follows. We first review the notion of
extended Nakajima monomials and the crystal structure given on the
set of such monomials. Also, we cite Young tableau expression of
crystal $\B(\infty)$ which play a crucial role in our work. We then
proceed to give monomial descriptions of $\B(\infty)$ and $\B(\la)$.
In the process of obtaining these results, we give new expressions
for the Kashiwara operators acting on a certain extended Nakajima
monomials, more appropriate for the situation in hand.

One strong point of the present work is that since our description
has a natural correspondence with the Young tableau description, the
Young tableau theory developed for $\B(\infty)$ and $\B(\la)$ can
easily be transferred to that on the monomial set.

In closing the introduction, we remark that recently the present
work has been extended to monomial descriptions of $\B(\infty)$ for
all other classical finite types and $G_2$ type by the
author~\cite{lee2}.

\vspace{2mm} \noindent \emph{Acknowledgments.} \ The author would
like to thank Professors Seok-Jin Kang, Satoshi Naito, and Yoshihisa
Saito for helpful information and valuable suggestions.

%%%%%%%%%%%%%%%%%%%%%%%%%%%%%%%%%%%%%%%%%%%%%%%%%%%%%%%%%%%%%%%%%%%%%%%%%%
\section{Extended Nakajima monomials and Young tableaux}%
\label{Nakajima monomials and crystals}

In this section, we introduce notation and cite facts that are crucial
for our work. Please refer to the references cited in the introduction
or books on quantum groups~\cite{MR1881971,MR1359532} for the basic
concepts on quantum groups and crystal bases.

Let us first fix the basic notation.
\begin{itemize}
\item $I=\{1,\dots,n\}$ : index set.
\item $A=(a_{ij})_{i,j\in I}$
      : Cartan matrix of type $\an$.
\item $P^{\vee} =\bigoplus_{i\in I} \Z h_i$ : dual weight lattice.
\item $P=\{\la \in \mathfrak h^* \mid \la(P^{\vee}) \subset
      \Z \}$ : weight lattice, where $\mathfrak h = \Q \ot_{\Z} P^{\vee}$.
\item $P^{+}=\{ \la \in P\vert \la(h_i) \ge 0 \ \text{for all} \
      i\in I\}$ : the set of dominant integral weights. \item
      $\Pi^{\vee} = \{h_i \mid i \in I\}$ : the set of simple coroots.
\item $\Pi = \{ \alpha_i \mid i \in I \}$ : the set of simple
      roots .
\item $\uq(\an)$ : quantum group associated with the Cartan datum
      $(A,\Pi,\Pi^\vee, P, P^\vee)$.
\item $\uq^-(\an)$ : subalgebra of $\uq(\an)$
      generated by $f_i$ ($i\in I$).
\item $\hwc(\la)$ : irreducible highest weight crystal of highest
      weight $\la$.
\item $\B(\infty)$ : crystal base of $\uq^-(\an)$.
\end{itemize}
Throughout this paper, a $U_q(\an)$-crystal will refer to a
(abstract) crystal associated with the Cartan datum
$(A,\Pi,\Pi^\vee, P, P^\vee)$. The crystal base $\B(\infty)$ of
$U_q^-(\an)$ is a $U_q(\an)$-crystal.

%%%%%%%%%%%%%%%%%%%%%%%%%%%%%%%%%%%%%%%%%%%%
\subsection{Nakajima monomials}

We now recall the set of monomials discovered by Nakajima and its
crystal structure and also recall their extension introduced
in~\cite{lee}. Both of these sets were defined for all symmetrizable
Kac-Moody algebras, but we shall restrict ourselves to the $\an$
case in this paper. Our exposition of the crystal structure on
Nakajima monomials follows that of Kashiwara~\cite{MR1988989}.

We denote by $\M$ the set of \emph{Nakajima monomials} in the
variables $Y_i(m)$ ($i\in I$, $m\in\Z$). That is
\begin{equation*}
\M =\left\{\prod_{(i,m)\in I\times\Z}\hspace{-4mm}
    {Y_i(m)}^{y_i(m)} \Big\vert
    \ y_i(m) \in\Z \ \textup{vanishes except at finitely many} \ (i,m)
    \right\}.
\end{equation*}
Fix any set of integers $c={(c_{ij})}_{i \neq j \in I}$ such that
$c_{ij}+c_{ji}=1$, and set
\begin{equation*}
A_i(m)=Y_i(m)Y_i(m+1)\prod_{j \neq i}
{Y_j(m+c_{ji})}^{\langle h_j,\alpha_i\rangle}.
\end{equation*}
The crystal structure on $\M$ is defined as follows. For every
monomial $M=\prod_{(i,m)\in I\times\Z}{Y_i(m)}^{y_i(m)}\in\M$, we
set
\begin{align*}
\wt(M) &= \sum_i(\sum_m y_i(m))\La_i,\\ \vphi_i(M) &=
\textup{max}\Big\{ \sum_{k\le m} y_i(k)
            \,\big\vert \ m\in\Z \Big\},\\
\veps_i(M) &= \textup{max}\Big\{ -\sum_{k>m} y_i(k)
            \,\big\vert \ m\in\Z \Big\}.
\end{align*}
We define
\begin{align}
\fit(M) &=
\begin{cases}
0 &\textup{if} \ \vphi_i(M)=0,\\
A_i(m_f)^{-1}M &\textup{if} \ \vphi_i(M)>0,
\end{cases}\label{f1}\\
\eit(M) &=
\begin{cases}
0 &\textup{if} \ \veps_i(M)=0,\\
A_i(m_e)M &\textup{if} \ \veps_i(M)>0.\label{e1}
\end{cases}
\end{align}
Here
\begin{align*}
m_f =\textup{min}\Big\{m \big\vert\  \vphi_i(M)=\sum_{k\le m}
y_i(k) \Big\},\quad m_e =\textup{max}\Big\{m \big\vert\
\veps_i(M)=-\sum_{k>m} y_i(k) \Big\}.
\end{align*}
These Kashiwara operators, together with the maps $\vphi_i$,
$\veps_i$ $(i\in I)$, $\wt$, define a crystal structure on the set
$\M$~\cite{MR1988989}. We denote by $\M_{c}$ the set $\M$ subject to
the crystal structure depending on the set $c$, as given above.

The following is a realization theorem for irreducible highest
weight crystal given by Kashiwara and Nakajima.

\begin{theorem}\textup{(}\cite{MR1988989}\textup{)}\label{nakakash}
For a maximal vector $M\in\M_{c}$, the connected component of
$\M_{c}$ containing $M$ is isomorphic to $\B(\wt(M))$.
\end{theorem}

\emph{Extended Nakajima monomials} and the crystal structure on
the set of such elements was introduced in~\cite{lee}.

Let $\Me$ be a certain set of formal monomials in the variables
$Y_i(m)^{(0,1)}$ and $Y_i(m)^{(1,0)}$ ($i\in I$, $m\in\Z$) given by
\begin{equation}
\Me =
\left\{\prod_{(i,m)\in I\times\Z}
        {Y_i(m)}^{y_i(m)}
        \Big\vert
\begin{aligned}
 & y_i(m)\!=\big(y^0_i(m),y^1_i(m)\big) \in\Z\times\Z\\
 & \textup{vanishes except at finitely many} \ (i,m)
\end{aligned}
\right\}.
\end{equation}
The product of monomials $Y_i(m)^{(u,v)}$ and $Y_i(m)^{(u',v')}$ are
set to $Y_i(m)^{(u+u',v+v')}$, for $(u,v)$, $(u',v')$
$\in\Z\times\Z$. We give the lexicographic order to the set
$\Z\times\Z$ of variable exponents.

Fix any set of integers $c={(c_{ij})}_{i \neq j \in I}$ such that
$c_{ij}+c_{ji}=1$, and set
\begin{equation}
A_i(m)^{\pm 1}={Y_i(m)}^{(0,\pm 1)}{Y_i(m+1)}^{(0,\pm 1)} \prod_{j
\neq i} {Y_j(m+c_{ji})}^{(0,\pm\langle h_j,\alpha_i\rangle)}.
\end{equation}

The crystal structure on $\Me$ is defined as follows. For each
monomial $M=\prod_{(i,m)\in I\times\Z}{Y_i(m)}^{y_i(m)}\in\Me$, we
set
\begin{align}
\wtt(M) &= \sum_i\big(\sum_m y_i(m)\big)\La_i
         = \sum_i\big(\sum_m (y^0_i(m),y^1_i(m))\big)\La_i,\\
\vphit_i(M) &= \textup{max}\Big\{ \sum_{k\le m} y_i(k)
            \,\Big\vert \ m\in\Z \Big\},\\
\vepst_i(M) &= \textup{max}\Big\{ -\sum_{k>m} y_i(k)
            \,\Big\vert \ m\in\Z \Big\}.
\end{align}
Notice that the coefficients of $\wtt(M)$ are pairs of integers.
In this setting, we have $\vphit_i(M)\geq (0,0)$, $\vepst_i(M)\geq
(0,0)$, and $\wtt(M)=\sum_i
\big(\vphit_i(M)-\vepst_i(M)\big)\La_i$. Set
\begin{align}
\wt(M) &= \sum_i\big(\sum_m y^1_i(m)\big)\La_i,\label{structure1}\\
\vphi_i(M) &= \sum_{k\le m} y^1_i(k) \quad \textup{where}\
              \vphit_i(M)= \sum_{k\le m} \big(y^0_i(k),y^1_i(k)\big),
              \label{structure2}\\
\veps_i(M) &= -\sum_{k>m} y^1_i(k) \quad \textup{where}\
              \vepst_i(M)= -\sum_{k>m} \big(y^0_i(k),y^1_i(k)\big).
              \label{structure3}
\end{align}
Then we trivially have $\wt(M)=\sum_i
\big(\vphi_i(M)-\veps_i(M)\big)\La_i$.
 From the above definition, ${Y_i(m)}^{(0,1)}$ has the weight $\La_i$,
and so $A_i(m)$ has the weight $\alpha_i$. We define the action of
Kashiwara operators by
\begin{align}
\fit(M) &=
\begin{cases}\label{actions}
0 &\textup{if} \ \vphit_i(M)=(0,0),\\
A_i(m_f)^{-1}M &\textup{if} \ \vphit_i(M)>(0,0),
\end{cases}\\
\eit(M) &=
\begin{cases}\label{actions2}
0 &\textup{if} \ \vepst_i(M)=(0,0),\\
A_i(m_e)M &\textup{if} \ \vepst_i(M)>(0,0).
\end{cases}
\end{align}
Here,
\begin{align}
m_f &=\textup{min}\Big\{m \big\vert\ \vphit_i(M)=
      \sum_{k\le m} y_i(k) \Big\},\label{emf}\\
m_e &=\textup{max}\Big\{m \big\vert\ \vepst_i(M)=
      -\sum_{k>m} y_i(k) \Big\}.\label{eme}
\end{align}
Note that $y_i(m_f)>(0,0)$, $y_i(m_f+1)\le (0,0)$,
$y_i(m_e+1)<(0,0)$, and $y_i(m_e)\ge (0,0)$.

For any fixed set of integers $c={(c_{ij})}_{i \neq j \in I}$ such
that $c_{ij}+c_{ji}=1$, the Kashiwara operators defined
in~\eqref{actions} and~\eqref{actions2}, together with the maps
$\vphi_i$, $\veps_i$ $(i\in I)$, and $\wt$ of~\eqref{structure1}
to~\eqref{structure3}, define a crystal structure on the set
$\Me$~\cite{lee}. We refer to an element of the set $\Me$ as an
\emph{extended Nakajima monomial} and denote by $\Mec$ the set $\Me$
subject to the crystal structure depending on the set $c$, as given
above.

\begin{remark}\label{naka}
Now, we may give many different crystal structures to the set of
extended Nakajima monomials through the choice of the set $c$. For
special linear Lie algebras, all the different crystals induced from
the set of extended Nakajima monomials through different choices of
the set $c$, are isomorphic (see ~\cite{MR1988989} or Proposition
3.2 of~\cite{lee}).

Consider the set $\dot\Mec$ of monomials $\prod_{(i,m)\in I\times\Z}
{Y_i(m)}^{(y_i^0(m),y_i^1(m))}\in\Mec$ with $y^0_i(m)=0$ for all
$(i,m)$. The set $\dot\Mec$ is a subcrystal of $\Mec$. It is exactly
the Nakajima monomial set $\M_c$ introduced in~\cite{MR1988989} if
we identify ${Y_i(m)}^{(0,y^1_i(m))}\in\dot\Mec$ with
${Y_i(m)}^{y^1_i(m)}\in\M_c$. The crystal structure on $\M_c$ is
compatible with that on $\Mec$ under this identification. We would
like to mention that $\M_c$ can be treated as a subcrystal of
$\Mec$. Viewing $\M_c$ as a subcrystal of $\Mec$, the monomial
theory developed for irreducible highest weight crystal can easily
be transferred to that on the extended monomial set.
\end{remark}

We obtain the following statement from Theorem~\ref{nakakash}.

\begin{cor}\textup{(}\cite{lee}\textup{)}\label{mla}
For a maximal vector $M\in\dot\Mec$, the connected component of
$\dot\Mec$ containing $M$ is crystal isomorphic to $\B(\wt(M))$.
\end{cor}

In the final section, we will give a concrete listing of elements
belonging to the connected component containing a maximal vector
$M\in\dot\Mec$ mentioned in the above corollary, for the case of
$\an$.

Unless there is possibility of confusion, we shall omit $c$ and use
the notations $\M$, $\dot\Me$, and $\Me$ instead of $\M_{c}$,
$\dot\Mec$, and $\Mec$, respectively.

%%%%%%%%%%%%%%%%%%%%%%%%%%%%%%%%%%%%%%%%%%%%
\subsection{Young tableaux}

In this section, we recall a Young tableau description of the
crystals $\B(\infty)$ and $\B(\la)$ for type $A_n$, that are
crucial for our work. Using these, in the remaining sections, we
shall show that the sets of monomials, satisfying some appropriate
conditions, give new descriptions of $\B(\infty)$ and $\B(\la)$.

Based on result of the paper~\cite{MR1273277}, we shall identify
elements of the highest weight crystal $\B(\la)$ with semistandard
tableaux of $\la$-shape, for the $A_n$ case. Since this work is a
well known result, we refer readers to the original papers and
shall not repeat the complicated definitions here.

\begin{definition}\hfill
\begin{enumerate}
\item A semistandard tableau $T$ of shape $\la\in P^+$,
      equivalently, an element of an irreducible highest weight
      crystal $\B(\la)$ for the $A_n$ type, is
      \emph{large} if it consists of $n$ non-empty rows, and if for each
      $1 \leq i \leq n$, the number of $i$-boxes in the $i$-th row is
      strictly greater than the number of all boxes in the $(i+1)$-th
      row. In particular the $n$-th row of $T$ contains at least one
      $n$-box.
\item A large tableau $T$ is \emph{marginally} large
      if for $1\leq i \leq n$, the number of $i$-boxes in the $i$-th row
      of $T$ is greater than the number of all boxes in the $(i+1)$-th
      row by exactly one. In particular, the $n$-th row of $T$ should
      contain one $n$-box.
\end{enumerate}
\end{definition}

In Figure~\ref{tbl:1}, for $A_3$ type, we give examples of
semistandard tableaux. The one on the left is large, the one in
the middle is marginally large, and the one on the right is not
large.
\begin{figure}
\centering

\begin{tabular}{rlll}
 \raisebox{-0.4\height}{
 \begin{texdraw}%
 \drawdim in
 \arrowheadsize l:0.065 w:0.03
 \arrowheadtype t:F
 \fontsize{6}{6}\selectfont
 \textref h:C v:C
 \drawdim em
 \setunitscale 1.4
 \move(0 0)
\bsegment \move(-6 3)\rlvec(6 0) \move(-6 2)\rlvec(6 0)\rlvec(0 1)
\move(-6 1)\rlvec(2 0)\rlvec(0 2) \move(-6 0)\rlvec(1 0)\rlvec(0
3) \move(-6 0)\rlvec(0 3) \move(-3 2)\rlvec(0 1) \move(-2
2)\rlvec(0 1) \move(-1 2)\rlvec(0 1) \htext(-0.5 2.5){$4$}
\htext(-1.5 2.5){$3$} \htext(-2.5 2.5){$1$} \htext(-3.5 2.5){$1$}
\htext(-4.5 2.5){$1$} \htext(-5.5 2.5){$1$} \htext(-4.5 1.5){$2$}
\htext(-5.5 1.5){$2$} \htext(-5.5 0.5){$3$} \esegment
 \end{texdraw}%
 }
 &
 \raisebox{-0.4\height}{
 \begin{texdraw}%
 \drawdim in
 \arrowheadsize l:0.065 w:0.03
 \arrowheadtype t:F
 \fontsize{6}{6}\selectfont
 \textref h:C v:C
 \drawdim em
 \setunitscale 1.4
 \move(-5 3)\rlvec(5 0) \move(-5 2)\rlvec(5 0)\rlvec(0 1) \move(-5
 1)\rlvec(3 0)\rlvec(0 2) \move(-5 0)\rlvec(1 0)\rlvec(0 3)
 \move(-5 0)\rlvec(0 3) \move(-3 1)\rlvec(0 2) \move(-1 2)\rlvec(0
 1) \move(-4 0)\rlvec(1 0)\rlvec(0 1)
 \htext(-0.5 2.5){$4$} \htext(-1.5 2.5){$1$} \htext(-2.5 2.5){$1$}
 \htext(-3.5 2.5){$1$} \htext(-4.5 2.5){$1$} \htext(-4.5 1.5){$2$}
 \htext(-3.5 1.5){$2$} \htext(-2.5 1.5){$2$} \htext(-4.5 0.5){$3$}
 \htext(-3.5 0.5){$4$}
 \end{texdraw}%
 }
 &
 \raisebox{-0.4\height}{
 \begin{texdraw}%
 \drawdim in
 \arrowheadsize l:0.065 w:0.03
 \arrowheadtype t:F
 \fontsize{6}{6}\selectfont
 \textref h:C v:C
 \drawdim em
 \setunitscale 1.4
 \move(0 2)\rlvec(4 0)
 \move(0 1)\rlvec(4 0)\rlvec(0 1)
 \move(0 0)\rlvec(2 0)\rlvec(0 2)
 \move(0 0)\rlvec(0 2)
 \move(1 0)\rlvec(0 2)
 \move(3 1)\rlvec(0 1)
 \htext(0.5 1.5){$1$}
 \htext(1.5 1.5){$1$}
 \htext(2.5 1.5){$2$}
 \htext(3.5 1.5){$3$}
 \htext(0.5 0.5){$3$}
 \htext(1.5 0.5){$4$}
 \end{texdraw}%
 }
\end{tabular}
\caption{Large (left), marginally large (middle), and non-large
(right) tableaux}\label{tbl:1}
\end{figure}

\begin{definition}
We denote by $\T(\infty)$ the set of all marginally large
tableaux. The marginally large tableau whose $i$-th row consists
only of $i$-boxes ($i\in I$) is denoted by $T_{\infty}$.
\end{definition}

\begin{example}
The set $\Tinf$ for case $A_3$, consists of all tableaux of the
following form. The unshaded part must exist, whereas the shaded
part is optional with variable size. {\allowdisplaybreaks
\begin{align*}
&T = \raisebox{-0.5\height}{\ %
\begin{texdraw}
 \fontsize{6}{6}\selectfont
 \textref h:C v:C
 \drawdim em
 \setunitscale 1.35
 \move(12 2)\rlvec(9 0)\rlvec(0 1)\rlvec(-9 0)
 \ifill f:0.8
 \move(5 1)\rlvec(6 0)\rlvec(0 2)\rlvec(-6 0)
 \ifill f:0.8
 \move(1 0)\rlvec(3 0)\rlvec(0 3)\rlvec(-3 0)
 \ifill f:0.8
\move(0 0) \bsegment
 \move(0 3)\rlvec(21 0)
 \move(0 2)\rlvec(21 0)\rlvec(0 1)
 \move(0 1)\rlvec(11 0)\rlvec(0 2)
 \move(0 0)\rlvec(4 0)\rlvec(0 2)
 \move(0 0)\rlvec(0 3)
 \move(1 0)\rlvec(0 1)
 \move(12 2)\rlvec(0 1)
 \move(5 1)\rlvec(0 1)
 \move(8 1)\rlvec(0 1)
 \move(15 2)\rlvec(0 1)
 \move(18 2)\rlvec(0 1)
\esegment
 \htext(0.5 0.5){$3$} \htext(1.5 0.5){$4$} \htext(2.5
 0.5){$\cdots$} \htext(3.5 0.5){$4$}
 \htext(0.5 1.5){$2$} \htext(2 1.5){$\cdots$} \htext(3.5 1.5){$2$}
 \htext(4.5 1.5){$2$} \htext(5.5 1.5){$3$} \htext(6.5 1.5){$\cdots$}
 \htext(7.5 1.5){$3$} \htext(8.5 1.5){$4$} \htext(9.5 1.5){$\cdots$}
 \htext(10.5 1.5){$4$}
 \htext(0.5 2.5){$1$} \htext(4.5 2.5){$1$} \htext(2.5 2.5){$\cdots$}
 \htext(10.5 2.5){$1$} \htext(7.5 2.5){$\cdots$}
 \htext(11.5 2.5){$1$} \htext(12.5 2.5){$2$} \htext(13.5 2.5){$\cdots$}
 \htext(14.5 2.5){$2$} \htext(15.5 2.5){$3$} \htext(16.5 2.5){$\cdots$}
 \htext(17.5 2.5){$3$} \htext(18.5 2.5){$4$} \htext(19.5 2.5){$\cdots$}
 \htext(20.5 2.5){$4$}
\end{texdraw}}\\
&T_{\infty} = \raisebox{-0.5\height}{\ %
\begin{texdraw}
 \fontsize{6}{6}\selectfont
 \textref h:C v:C
 \drawdim em
 \setunitscale 1.35
 \move(0 0)\lvec(3 0)\lvec(3 -1)\lvec(0 -1)
 \move(2 0)\lvec(2 -2)\lvec(0 -2)
 \move(1 0)\lvec(1 -3)\lvec(0 -3)
 \move(0 0)\lvec(0 -3)
 \htext(0.5 -0.5){$1$}
 \htext(1.5 -0.5){$1$}
 \htext(2.5 -0.5){$1$}
 \htext(0.5 -1.5){$2$}
 \htext(1.5 -1.5){$2$}
 \htext(0.5 -2.5){$3$}
\end{texdraw}}\,
\end{align*}}% end of allowdisplaybreaks
\end{example}

We recall the action of Kashiwara operators $\fit$, $\eit$ $(i\in
I)$ on marginally large tableaux $T\in\T(\infty)$.
\begin{enumerate}
\item We first read the boxes in the tableau $T$ through
      the \emph{far eastern reading} and write down the boxes
      in \emph{tensor product form.}
      That is, we read through each column from top to bottom starting
      from the rightmost column, continuing to the left, and lay down
      the read boxes from left to right in tensor product form.
      The following diagram gives an example.
\begin{equation*}
\qquad\ \raisebox{-0.4\height}{
\begin{texdraw}%
\fontsize{6}{6}\selectfont \textref h:C v:C \drawdim em
\setunitscale 1.4
\move(0 3)\rlvec(5 0) \move(0 2)\rlvec(5 0)\rlvec(0 1) \move(0
1)\rlvec(4 0)\rlvec(0 2) \move(0 0)\rlvec(2 0)\rlvec(0 2) \move(0
0)\rlvec(0 3) \move(1 0)\rlvec(0 3) \move(3 1)\rlvec(0 2) \move(2
1)\rlvec(0 2) \htext(0.5 2.5){$1$} \htext(1.5 2.5){$1$} \htext(2.5
2.5){$1$} \htext(3.5 2.5){$1$} \htext(4.5 2.5){$1$} \htext(0.5
1.5){$2$} \htext(1.5 1.5){$2$} \htext(2.5 1.5){$2$} \htext(3.5
1.5){$4$} \htext(0.5 0.5){$3$} \htext(1.5 0.5){$4$}
\end{texdraw}%
} = \raisebox{-0.3\height}{
\begin{texdraw}%
\fontsize{6}{6}\selectfont \textref h:C v:C \drawdim em
\setunitscale 1.4
\move(0 0)\rlvec(1 0)\rlvec(0 1)\rlvec(-1 0)\rlvec(0 -1)
\htext(0.5 0.5){$1$}
\end{texdraw}%
} \otimes
\raisebox{-0.3\height}{%
\begin{texdraw}%
\fontsize{6}{6}\selectfont \textref h:C v:C \drawdim em
\setunitscale 1.4
\move(0 0)\rlvec(1 0)\rlvec(0 1)\rlvec(-1 0)\rlvec(0 -1)
\htext(0.5 0.5){$1$}
\end{texdraw}%
} \otimes
\raisebox{-0.3\height}{%
\begin{texdraw}%
\fontsize{6}{6}\selectfont \textref h:C v:C \drawdim em
\setunitscale 1.4
\move(0 0)\rlvec(1 0)\rlvec(0 1)\rlvec(-1 0)\rlvec(0 -1)
\htext(0.5 0.5){$4$}
\end{texdraw}%
} \otimes
\raisebox{-0.3\height}{%
\begin{texdraw}%
\fontsize{6}{6}\selectfont \textref h:C v:C \drawdim em
\setunitscale 1.4
\move(0 0)\rlvec(1 0)\rlvec(0 1)\rlvec(-1 0)\rlvec(0 -1)
\htext(0.5 0.5){$1$}
\end{texdraw}%
} \otimes
\raisebox{-0.3\height}{%
\begin{texdraw}%
\fontsize{6}{6}\selectfont \textref h:C v:C \drawdim em
\setunitscale 1.4
\move(0 0)\rlvec(1 0)\rlvec(0 1)\rlvec(-1 0)\rlvec(0 -1)
\htext(0.5 0.5){$2$}
\end{texdraw}%
} \otimes
\raisebox{-0.3\height}{%
\begin{texdraw}%
\fontsize{6}{6}\selectfont \textref h:C v:C \drawdim em
\setunitscale 1.4
\move(0 0)\rlvec(1 0)\rlvec(0 1)\rlvec(-1 0)\rlvec(0 -1)
\htext(0.5 0.5){$1$}
\end{texdraw}%
} \otimes
\raisebox{-0.3\height}{%
\begin{texdraw}%
\fontsize{6}{6}\selectfont \textref h:C v:C \drawdim em
\setunitscale 1.4
\move(0 0)\rlvec(1 0)\rlvec(0 1)\rlvec(-1 0)\rlvec(0 -1)
\htext(0.5 0.5){$2$}
\end{texdraw}%
} \otimes
\raisebox{-0.3\height}{%
\begin{texdraw}%
\fontsize{6}{6}\selectfont \textref h:C v:C \drawdim em
\setunitscale 1.4
\move(0 0)\rlvec(1 0)\rlvec(0 1)\rlvec(-1 0)\rlvec(0 -1)
\htext(0.5 0.5){$4$}
\end{texdraw}%
} \otimes
\raisebox{-0.3\height}{%
\begin{texdraw}%
\fontsize{6}{6}\selectfont \textref h:C v:C \drawdim em
\setunitscale 1.4
\move(0 0)\rlvec(1 0)\rlvec(0 1)\rlvec(-1 0)\rlvec(0 -1)
\htext(0.5 0.5){$1$}
\end{texdraw}%
} \otimes
\raisebox{-0.3\height}{%
\begin{texdraw}%
\fontsize{6}{6}\selectfont \textref h:C v:C \drawdim em
\setunitscale 1.4
\move(0 0)\rlvec(1 0)\rlvec(0 1)\rlvec(-1 0)\rlvec(0 -1)
\htext(0.5 0.5){$2$}
\end{texdraw}%
} \otimes
\raisebox{-0.3\height}{%
\begin{texdraw}%
\fontsize{6}{6}\selectfont \textref h:C v:C \drawdim em
\setunitscale 1.4
\move(0 0)\rlvec(1 0)\rlvec(0 1)\rlvec(-1 0)\rlvec(0 -1)
\htext(0.5 0.5){$3$}
\end{texdraw}%
}
\end{equation*}
\item Under each tensor component $x$ of $T$,
      write down $\veps_i(x)$-many 1s followed by $\vphi_i(x)$-many 0s.
      Then, from the long sequence of mixed 0s and 1s,
      successively cancel out every occurrence of (0,1) pair
      until we arrive at a sequence of 1s followed by 0s,
      reading from left to right. This is called the
      $i$-signature of $T$.
\item Denote by $T'$, the tableau obtained from $T$,
      by replacing the box $x$ corresponding to the leftmost $0$
      in the $i$-signature of $T$ with the box $\fit x$.
      \begin{itemize}
      \item If $T'$ is a large tableau,
            it is automatically marginally large.
            We define $\fit T$ to be $T'$.
      \item If $T'$ is not large, then we define $\fit T$ to be the
            large tableau obtained by inserting one column
            consisting of $i$ rows to the left of the box $\fit$
            acted upon. The added column should have a
            $k$-box at the $k$-th row for $1\leq k \leq i$.
      \end{itemize}
\item Denote by $T'$, the tableau obtained from $T$,
      by replacing the box $x$
      corresponding to the rightmost $1$ in the $i$-signature of
      $T$ with the box $\eit x$.
      \begin{itemize}
      \item If $T'$ is a marginally large tableau,
            then we define $\eit T$ to be $T'$.
      \item If $T'$ is large but not marginally large,
            then we define $\eit T$ to be the
            large tableau obtained by removing the column
            containing the changed box.
            It will be of $i$ rows and have a $k$-box
            at the $k$-th row for $1\leq k \leq i$.
      \end{itemize}
\item If there is no $1$ in the $i$-signature of $T$,
      we define $\eit T=0$.
\end{enumerate}

\begin{remark}
The condition \emph{large} imposed on the tableau $T$ ensures that
its $i$-signature always contains $0$'s.
\end{remark}

Let $T$ be a tableau in $\T(\infty)$ with the $i$-th row, for each
$1\leq i \leq n$, consisting of $b^i_j$-many $j$\,s ($i< j\le
n+1$) and some number of $i$\,s. We define the maps $\wt:
\T(\infty) \rightarrow P$, $\vphi_i, \veps_i : \T(\infty)
\rightarrow \Z$ by setting
\begin{align*}
\wt(T)      &=-\sum_{j=1}^n\Big(\sum_{k=j+1}^{n+1}b_{k}^1+
                            \sum_{k=j+1}^{n+1}b_{k}^2+
                            \cdots +
                            \sum_{k=j+1}^{n+1}b_{k}^j\Big)\alpha_j,\\
\veps_i(T)  &=\text{the number of $1$s
             in the $i$-signature of $T$},\\
\vphi_i(T)  &=\veps_i(T)+\langle h_i,\wt(T)\rangle.
\end{align*}

\begin{theorem}\textup{(}\cite{HL}\textup{)}\label{HL}
The Kashiwara operators and various maps given above define a
crystal structure on $\Tinf$. The crystal $\T(\infty)$ is
isomorphic to $\B(\infty)$ as $U_q(A_n)$-crystal.
\end{theorem}

\begin{figure}
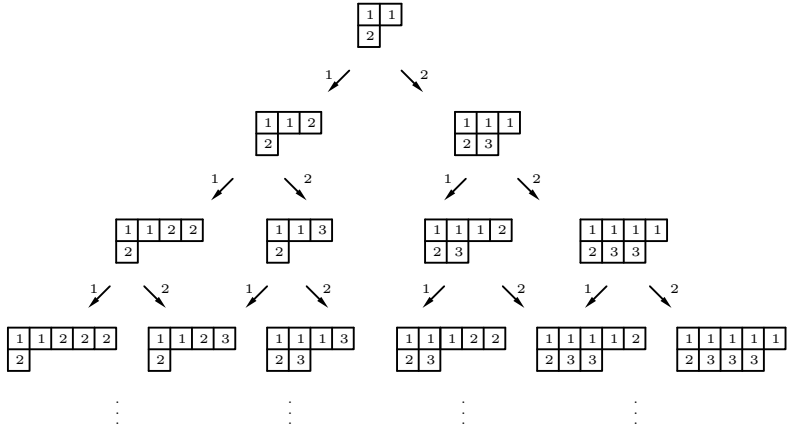

\centering
\begin{texdraw}%
\drawdim in \arrowheadsize l:0.065 w:0.03 \arrowheadtype t:F
\fontsize{5}{5}\selectfont \textref h:C v:C \drawdim em
\setunitscale 1.2
\move(-1.2 -1.1)\ravec(-1 -1)%%%%%%%%%%%%%%%
\move(1.2 -1.1)\ravec(1 -1)%%%%%%%%%%%%%%%%%
\move(-6.6 -6.1)\ravec(-1 -1)%%%%%%%%%%%%%%%
\move(-4.2 -6.1)\ravec(1 -1)%%%%%%%%%%%%%%%%
\move(4.2 -6.1)\ravec(-1 -1)%%%%%%%%%%%%%%%%
\move(6.6 -6.1)\ravec(1 -1)%%%%%%%%%%%%%%%%%
\move(-12.3 -11.1)\ravec(-1 -1)%%%%%%%%%%%%%
\move(-10.7 -11.1)\ravec(1 -1)%%%%%%%%%%%%%%
\move(-5 -11.1)\ravec(-1 -1)%%%%%%%%%%%%%%%%
\move(-3.2 -11.1)\ravec(1 -1)%%%%%%%%%%%%%%%
\move(3.2 -11.1)\ravec(-1 -1)%%%%%%%%%%%%%%%
\move(5.9 -11.1)\ravec(1 -1)%%%%%%%%%%%%%%%%
\move(10.7 -11.1)\ravec(-1 -1)%%%%%%%%%%%%%%
\move(12.7 -11.1)\ravec(1 -1)%%%%%%%%%%%%%%%
\htext(-12 -16.6){$\vdots$}%%%%%%%%%%%%%%%%%
\htext(4 -16.6){$\vdots$}%%%%%%%%%%%%%%%%%%%
\htext(-4 -16.6){$\vdots$}%%%%%%%%%%%%%%%%%%
\htext(12 -16.6){$\vdots$}%%%%%%%%%%%%%%%%%%
%%%%%%%%%%%%%%%%%%%%%%%%%%%%%%%%%%%
\move(1.2 0) \bsegment \move(-2 2)\rlvec(2 0) \move(-2 1)\rlvec(2
0)\rlvec(0 1) \move(-2 0)\rlvec(1 0)\rlvec(0 2) \move(-2
2)\rlvec(0 -2) \htext(-0.5 1.5){$1$} \htext(-1.5 1.5){$1$}
\htext(-1.5 0.5){$2$} \esegment
%%%%%%%%%%%%%%%%%%%%%%%%%%%%%%%%%%%
\move(-2.5 -5) \bsegment \move(-3 2)\rlvec(3 0) \move(-3
1)\rlvec(3 0)\rlvec(0 1) \move(-3 0)\rlvec(1 0)\rlvec(0 2)
\move(-3 2)\rlvec(0 -2) \move(-1 2)\rlvec(0 -1) \htext(-1.5
1.5){$1$} \htext(-2.5 1.5){$1$} \htext(-2.5 0.5){$2$} \htext(-0.5
1.5){$2$}\esegment
%%%%%%%%%%%%%%%%
\move(6.7 -5) \bsegment \move(-3 2)\rlvec(3 0) \move(-3 1)\rlvec(3
0)\rlvec(0 1) \move(-3 0)\rlvec(2 0)\rlvec(0 2) \move(-2
2)\rlvec(0 -2) \move(-3 2)\rlvec(0 -2) \htext(-0.5 1.5){$1$}
\htext(-1.5 1.5){$1$} \htext(-2.5 1.5){$1$} \htext(-2.5 0.5){$2$}
\htext(-1.5 0.5){$3$} \esegment
%%%%%%%%%%%%%%%%%%%%%%%%%%%%%%%%%%%%%%%%%%
\move(-8 -10) \bsegment \move(-4 2)\rlvec(4 0) \move(-4 1)\rlvec(4
0)\rlvec(0 1) \move(-4 0)\rlvec(1 0)\rlvec(0 2) \move(-4
2)\rlvec(0 -2) \move(-1 2)\rlvec(0 -1) \move(-2 2)\rlvec(0 -1)
\htext(-0.5 1.5){$2$} \htext(-1.5 1.5){$2$} \htext(-2.5 1.5){$1$}
\htext(-3.5 1.5){$1$} \htext(-3.5 0.5){$2$} \esegment
%%%%%%%%%%%%%%%%
\move(-2 -10) \bsegment \move(-3 2)\rlvec(3 0) \move(-3 1)\rlvec(3
0)\rlvec(0 1) \move(-3 0)\rlvec(1 0)\rlvec(0 2) \move(-3
2)\rlvec(0 -2) \move(-1 2)\rlvec(0 -1) \htext(-1.5 1.5){$1$}
\htext(-2.5 1.5){$1$} \htext(-2.5 0.5){$2$} \htext(-0.5 1.5){$3$}
\esegment
%%%%%%%%%%%%%%%%
\move(6.3 -10) \bsegment \move(-4 2)\rlvec(4 0) \move(-4
1)\rlvec(4 0)\rlvec(0 1) \move(-4 0)\rlvec(2 0)\rlvec(0 2)
\move(-4 2)\rlvec(0 -2) \move(-3 2)\rlvec(0 -2) \move(-1
2)\rlvec(0 -1) \htext(-0.5 1.5){$2$} \htext(-1.5 1.5){$1$}
\htext(-2.5 1.5){$1$} \htext(-3.5 1.5){$1$} \htext(-3.5 0.5){$2$}
\htext(-2.5 0.5){$3$} \esegment
%%%%%%%%%%%%%%%%
\move(13.5 -10) \bsegment \move(-4 2)\rlvec(4 0) \move(-4
1)\rlvec(4 0)\rlvec(0 1) \move(-4 0)\rlvec(3 0)\rlvec(0 2)
\move(-4 2)\rlvec(0 -2) \move(-2 2)\rlvec(0 -2) \move(-3
2)\rlvec(0 -2) \htext(-0.5 1.5){$1$} \htext(-1.5 1.5){$1$}
\htext(-2.5 1.5){$1$} \htext(-3.5 1.5){$1$} \htext(-3.5 0.5){$2$}
\htext(-2.5 0.5){$3$} \htext(-1.5 0.5){$3$} \esegment
%%%%%%%%%%%%%%%%%%%%%%%%%%%%%%%%%%%%%%%%%%%%%%%%%%%%
\move(-13 -15) \bsegment \move(0 1)\rlvec(1 0)\rlvec(0 1)\rlvec(-1
0) \move(-4 2)\rlvec(4 0) \move(-4 1)\rlvec(4 0)\rlvec(0 1)
\move(-4 0)\rlvec(1 0)\rlvec(0 2) \move(-4 2)\rlvec(0 -2) \move(-1
2)\rlvec(0 -1) \move(-2 2)\rlvec(0 -1) \htext(-0.5 1.5){$2$}
\htext(0.5 1.5){$2$} \htext(-1.5 1.5){$2$} \htext(-2.5 1.5){$1$}
\htext(-3.5 1.5){$1$} \htext(-3.5 0.5){$2$} \esegment
%%%%%%%%%%%%%%%%
\move(-6.5 -15) \bsegment \move(-4 2)\rlvec(4 0) \move(-4
1)\rlvec(4 0)\rlvec(0 1) \move(-4 0)\rlvec(1 0)\rlvec(0 2)
\move(-4 2)\rlvec(0 -2) \move(-1 2)\rlvec(0 -1) \move(-2
2)\rlvec(0 -1) \htext(-0.5 1.5){$3$} \htext(-1.5 1.5){$2$}
\htext(-2.5 1.5){$1$} \htext(-3.5 1.5){$1$} \htext(-3.5 0.5){$2$}
\esegment
%%%%%%%%%%%%%%%%
\move(-1 -15) \bsegment \move(-4 2)\rlvec(4 0) \move(-4 1)\rlvec(4
0)\rlvec(0 1) \move(-4 0)\rlvec(2 0)\rlvec(0 2) \move(-4
2)\rlvec(0 -2) \move(-3 2)\rlvec(0 -2) \move(-1 2)\rlvec(0 -1)
\htext(-0.5 1.5){$3$} \htext(-1.5 1.5){$1$} \htext(-2.5 1.5){$1$}
\htext(-3.5 1.5){$1$} \htext(-3.5 0.5){$2$} \htext(-2.5 0.5){$3$}
\esegment
%%%%%%%%%%%%%%%%
\move(5 -15) \bsegment \move(-4 2)\rlvec(5 0) \move(-4 1)\rlvec(5
0)\rlvec(0 1) \move(-4 0)\rlvec(2 0)\rlvec(0 2) \move(-3
2)\rlvec(0 -2) \move(-4 2)\rlvec(0 -2) \move(-1 2)\rlvec(0 -1)
\move(-2 2)\rlvec(0 -1) \move(0 2)\rlvec(0 -1) \htext(0.5
1.5){$2$} \htext(-0.5 1.5){$2$} \htext(-1.5 1.5){$1$} \htext(-2.5
1.5){$1$} \htext(-3.5 1.5){$1$} \htext(-3.5 0.5){$2$} \htext(-2.5
0.5){$3$} \esegment
%%%%%%%%%%%%%%%%
\move(11.5 -15) \bsegment \move(-4 2)\rlvec(5 0) \move(-4
1)\rlvec(5 0)\rlvec(0 1) \move(-4 0)\rlvec(3 0)\rlvec(0 2)
\move(-4 2)\rlvec(0 -2) \move(-3 2)\rlvec(0 -2) \move(-2
2)\rlvec(0 -2) \move(0 2)\rlvec(0 -1) \htext(0.5 1.5){$2$}
\htext(-0.5 1.5){$1$} \htext(-1.5 1.5){$1$} \htext(-2.5 1.5){$1$}
\htext(-3.5 1.5){$1$} \htext(-3.5 0.5){$2$} \htext(-2.5 0.5){$3$}
\htext(-1.5 0.5){$3$} \esegment
%%%%%%%%%%%%%%%%
\move(18 -15) \bsegment \move(-4 2)\rlvec(5 0) \move(-4 1)\rlvec(5
0)\rlvec(0 1) \move(-4 0)\rlvec(3 0)\rlvec(0 2) \move(-4
2)\rlvec(0 -2) \move(-3 2)\rlvec(0 -2) \move(-2 2)\rlvec(0 -2)
\move(0 2)\rlvec(0 -1) \move(-1 0)\rlvec(1 0)\rlvec(0 1)
\htext(0.5 1.5){$1$} \htext(-0.5 1.5){$1$} \htext(-1.5 1.5){$1$}
\htext(-2.5 1.5){$1$} \htext(-3.5 1.5){$1$} \htext(-3.5 0.5){$2$}
\htext(-2.5 0.5){$3$} \htext(-0.5 0.5){$3$} \htext(-1.5 0.5){$3$}
\esegment
%%%%%%%%%%%%%%%%%%%%%%%%%%%%%%%%%%%%%%%%%%%%%%%%%%%%
\move(0 0) \bsegment \htext(-2.2 -1.3){$1$} \htext(2.2 -1.3){$2$}
\esegment \move(0 -2.7) \bsegment \htext(-7.5 -3.4){$1$}
\htext(-3.2 -3.4){$2$} \htext(3.3 -3.4){$1$} \htext(7.4 -3.4){$2$}
\esegment \move(0 -6.2) \bsegment \htext(-9.8 -5){$2$}
\htext(-13.1 -5){$1$} \htext(-5.9 -5){$1$} \htext(-2.3 -5){$2$}
\htext(2.3 -5){$1$} \htext(6.7 -5){$2$} \htext(9.8 -5){$1$}
\htext(13.8 -5){$2$} \esegment \move(0 1) \move(0 -14)
%\drawbb
\end{texdraw}%
\caption{Part of the crystal $\T(\infty)$ for type
$A_2$}\label{fig:1}
\end{figure}

In the figure~\ref{fig:1}, we illustrate the top part of the
crystal $\T(\infty)$ for type $A_2$.

%%%%%%%%%%%%%%%%%%%%%%%%%%%%%%%%%%%%%%%%%%%%%%%%%%%%%%%%%%%%%%%%%%%%%%%%%%
\section{Monomial description of crystal $\B(\infty)$}%
\label{Extended monomial description of crystals
B(infty) and B(la)}

We give an explicit description of the crystal $\B(\infty)$ for
$\an$-type, in terms of extended Nakajima monomials. We first
present a candidate monomial set, show this set to be a crystal, and
give a crystal isomorphism of this with another description of
$\B(\infty)$.

We take the set $c=(c_{ij})_{i\neq j}$ to be
\begin{equation}\label{setc}
c_{ij}=
\begin{cases}
      0 & \quad \textup{ if } i>j,\\
      1 & \quad \textup{ if } i<j.
\end{cases}
\end{equation}
Then for $i\in I$ and $m\in\Z$, we have
\begin{equation*}
A_i(m)^{\pm 1}={Y_i(m)}^{(0,\pm 1)}{Y_i(m+1)}^{(0,\pm 1)}
        {Y_{i-1}(m+1)}^{(0,\mp 1)}{Y_{i+1}(m)}^{(0,\mp 1)}.
\end{equation*}
Here, we are setting $Y_0(k)^{(0,\pm 1)}=Y_{n+1}(k)^{(0,\pm 1)}=1$.

As was stated in Remark~\ref{naka}, the theory to be developed on
the fixed $\Mec$ can similarly be developed on all the other
crystals induced from the extended Nakajima monomial set through
other choices of the set $c$, under the isomorphism given in
Proposition 3.2 (2) of~\cite{lee}.
 From now on, we shall omit $c$ and use the notation $\Me$
instead of $\Mec$, since we already fixed the set $c$.

The set we define below is originally obtained by applying Kashiwara
operators $\fit$ ($i\in I$) iteratively, starting from the maximal
element $\prod_{i\in I} {Y_i(-i)}^{(1,0)}\in\Me$ of
$\wtt(\prod_{i\in I} {Y_i(-i)}^{(1,0)})=\sum_i(1,0)\La_i$. This
choice of starting monomial will allow us to relate monomials of the
set defined below to tableaux in $\T(\infty)$ naturally.

\begin{definition}
Consider elements of $\Me$ having the form
\begin{equation}\label{extend}
\begin{aligned}
      M=& \prod_{i\in I} \Big({Y_i(-i)}^{(1,a_i^i)}
          \prod_{m=0}^{i-1}{Y_i(-m)}^{(0,a_i^m)}\Big)\\
      =        &{Y_1(-1)}^{(1,a_1^1)}{Y_1(0)}^{(0,a_1^0)}\\
        &\cdot {Y_2(-2)}^{(1,a_2^2)}{Y_2(-1)}^{(0,a_2^1)}
                {Y_2(0)}^{(0,a_2^0)}\\
        & \quad\cdots\\
        &\cdot {Y_n(-n)}^{(1,a_n^n)}{Y_n(-n\!+\!1)}^{(0,a_n^{n\!-\!1})}
        \cdots {Y_n(-1)}^{(0,a_n^1)}{Y_n(0)}^{(0,a_n^0)}
\end{aligned}
\end{equation}
with the conditions
\begin{enumerate}
\item $\sum_{j=0}^k a_{i+j}^j\le 0$
      for each $0\le k\le n\!-\!1$, $1\le i\le n\!-\!k$,
\item $\sum_{i=1}^{n\!-\!k}(\sum_{j=0}^k
                  a_{i+j}^j)=\sum_{i=k+1}^n a_i^i$
      for $0\le k\le n\!-\!1$.
\end{enumerate}
Specifically, in case of $a_i^j=0$ for all $i,j$, we have
\begin{equation}\label{hi9}
      M=\prod_{i\in I} {Y_i(-i)}^{(1,0)}
       ={Y_1(-1)}^{(1,0)}{Y_2(-2)}^{(1,0)}\cdots {Y_n(-n)}^{(1,0)}.
\end{equation}
We denote by $\M(\infty)$ the set of all monomials of the
form~(\ref{extend}) and by $M_{\infty}$ the monomial
of~\eqref{hi9}.
\end{definition}

\begin{figure}
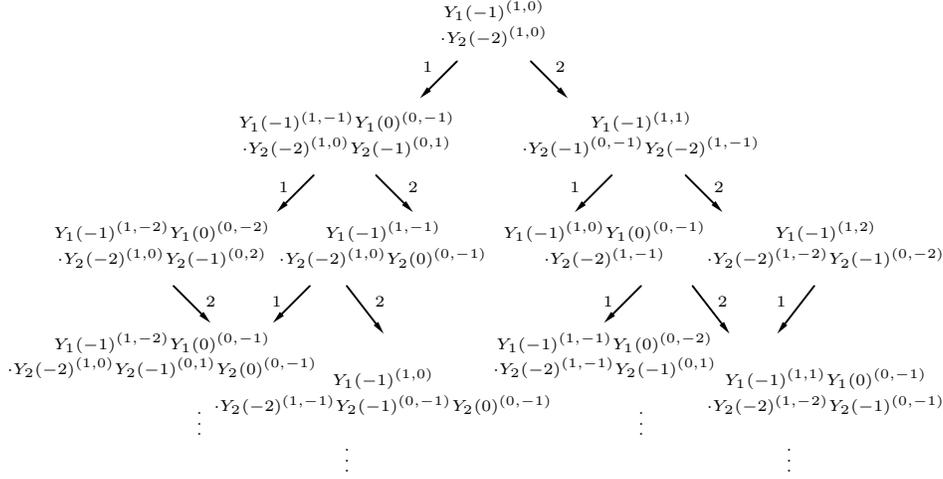

\centering
\begin{texdraw}%
\drawdim in \arrowheadsize l:0.065 w:0.03 \arrowheadtype t:F
\fontsize{6}{6}\selectfont \textref h:C v:C \drawdim em
\setunitscale 1.9
\move(-1 -1.4)\ravec(-1 -1)%%%%%%%%%%%%%%%
\move(1 -1.4)\ravec(1 -1)%%%%%%%%%%%%%%%%%
\move(-3.2 -4.5)\ravec(1 -1)%%%%%%%%%%%%%%
\move(-4.9 -4.5)\ravec(-1 -1)%%%%%%%%%%%%%
\move(3.2 -4.5)\ravec(-1 -1)%%%%%%%%%%%%%%
\move(5.2 -4.5)\ravec(1 -1)%%%%%%%%%%%%%%%
\move(-4 -7.5)\ravec(1 -1.3)%%%%%%%%%%%%%%
\move(-5 -7.5)\ravec(-1 -1)%%%%%%%%%%%%%%%
\move(8.7 -7.5)\ravec(-1 -1.3)%%%%%%%%%%%%
\move(-8.7 -7.5)\ravec(1 -1)%%%%%%%%%%%%%%
\move(4 -7.5)\ravec(-1 -1)%%%%%%%%%%%%%%%%
\move(5.4 -7.5)\ravec(1 -1.3)%%%%%%%%%%%%%
\htext(-8 -11){$\vdots$}%%%%%%%%%%%%%%%%
\htext(4 -11){$\vdots$}%%%%%%%%%%%%%%%%%
\htext(-4 -12){$\vdots$}%%%%%%%%%%%%%%%%
\htext(8 -12){$\vdots$}%%%%%%%%%%%%%%%%%

\htext(0 0){${Y_1(-1)}^{(1,0)}$} \htext(0 -0.7){$\cdot
{Y_2(-2)}^{(1,0)}$}

\htext(-4 -3){${Y_1(-1)}^{(1,-1)}{Y_1(0)}^{(0,-1)}$} \htext(-4
-3.7){$\cdot {Y_2(-2)}^{(1,0)}{Y_2(-1)}^{(0,1)}$}

\htext(4 -3){${Y_1(-1)}^{(1,1)}$} \htext(4 -3.7){$\cdot
{Y_2(-1)}^{(0,-1)}{Y_2(-2)}^{(1,-1)}$}

\htext(-9 -6){${Y_1(-1)}^{(1,-2)}{Y_1(0)}^{(0,-2)}$} \htext(-9
-6.7){$\cdot {Y_2(-2)}^{(1,0)}{Y_2(-1)}^{(0,2)}$}

\htext(-3 -6){${Y_1(-1)}^{(1,-1)}$} \htext(-3 -6.7){$\cdot
{Y_2(-2)}^{(1,0)}{Y_2(0)}^{(0,-1)}$}

\htext(3 -6){${Y_1(-1)}^{(1,0)}{Y_1(0)}^{(0,-1)}$} \htext(3
-6.7){$\cdot {Y_2(-2)}^{(1,-1)}$}

\htext(9 -6){${Y_1(-1)}^{(1,2)}$} \htext(9 -6.7){$\cdot
{Y_2(-2)}^{(1,-2)}{Y_2(-1)}^{(0,-2)}$}

\htext(-9 -9){${Y_1(-1)}^{(1,-2)}{Y_1(0)}^{(0,-1)}$} \htext(-9
-9.7){$\cdot {Y_2(-2)}^{(1,0)}{Y_2(-1)}^{(0,1)}
                       {Y_2(0)}^{(0,-1)}$}

\htext(-3 -10){${Y_1(-1)}^{(1,0)}$} \htext(-3 -10.7){$\cdot
{Y_2(-2)}^{(1,-1)}{Y_2(-1)}^{(0,-1)}
                       {Y_2(0)}^{(0,-1)}$}

\htext(3 -9){${Y_1(-1)}^{(1,-1)}{Y_1(0)}^{(0,-2)}$} \htext(3
-9.7){$\cdot {Y_2(-2)}^{(1,-1)}{Y_2(-1)}^{(0,1)}$}

\htext(9 -10){${Y_1(-1)}^{(1,1)}{Y_1(0)}^{(0,-1)}$} \htext(9
-10.7){$\cdot {Y_2(-2)}^{(1,-2)}{Y_2(-1)}^{(0,-1)}$}
%%%%%%%%%%%%%%
\move(0 0) \bsegment \htext(-1.8 -1.5){$1$} \htext(1.8 -1.5){$2$}
\esegment \move(0 -2.7) \bsegment \htext(-5.7 -2.1){$1$}
\htext(-2.2 -2.1){$2$} \htext(2.2 -2.1){$1$} \htext(6.1 -2.1){$2$}
\esegment \move(0 -6.2) \bsegment \htext(-7.7 -1.7){$2$}
\htext(-5.9 -1.7){$1$} \htext(-3.1 -1.7){$2$} \htext(3.1
-1.7){$1$} \htext(6.2 -1.7){$2$} \htext(7.8 -1.7){$1$} \esegment
%\move(0 -13.3)
%\drawbb
\end{texdraw}%
\caption{Part of the monomial set $\M(\infty)$ for $A_2$
type}\label{fig:2}
\end{figure}

Actually, as we will become apparent later, this set $\M(\infty)$ is
closed and connected under Kashiwara operators~(\ref{actions})
and~(\ref{actions2}) on $\Me$. Figure~\ref{fig:2} is the top part of
monomial set $\M(\infty)$.

We now introduce new expressions for elements of $\M(\infty)$.
First, we introduce the following notation.

\begin{definition}\label{xxx}
For $u,v,m\in\Z$, and $1\le i\le n\!+\!1$, we use the notation
\begin{equation*}%\label{order}
{X_i(m)}^{(u,v)} ={Y_i(m)}^{(u,v)}{Y_{i-1}(m+1)}^{(-u,-v)}.
\end{equation*}
Here, we set $Y_0(k)^{(u,v)}=Y_{n+1}(k)^{(u,v)}=1$.
\end{definition}

\begin{remark}
 From the above notation, we obtain
\begin{equation}\label{1}\begin{aligned}
Y_i(m)^{(u,v)}&=X_1(m+i-1)^{(u,v)}
                X_2(m+i-2)^{(u,v)}
                \cdots
                X_i(m)^{(u,v)}\\
              &=X_{n+1}(m+i-(n+1))^{(-u,-v)}
                X_n(m+i-n)^{(-u,-v)}\\
              &\quad\cdots
                X_{i+1}(m-1)^{(-u,-v)},
\end{aligned}\end{equation}
for each $i\in I$. And so, we may write
\begin{equation*}
A_i(m) = X_i(m)^{(0,1)}X_{i+1}(m)^{(0,-1)}.
\end{equation*}
This is very useful when computing Kashiwara action on monomials
written in terms of $X_i(m)^{(u,v)}$.
\end{remark}

\begin{prop}
Consider elements of $\Me$ having the form
\begin{equation}~\label{equ:jin8}
\begin{aligned}
      M=& \prod_{i\in I}\Big(X_i(-i)^{(n-i+1,-\sum_{k=i+1}^{n\!+\!1}b_k^i)}
          \prod_{k=i+1}^{n+1}X_k(-i)^{(0,b_k^i)}\Big)\\
       =& X_1(-1)^{(n,-\sum_{k=2}^{n\!+\!1}b_k^1)}
          X_2(-1)^{(0,b_2^1)}
          \cdots
          X_n(-1)^{(0,b_n^1)}
          X_{n+1}(-1)^{(0,b_{n\!+\!1}^1)}\\
          \cdot\,
        & X_2(-2)^{(n\!-\!1,-\sum_{k=3}^{n\!+\!1}b_k^2)}
          X_3(-2)^{(0,b_3^2)}
          \cdots
          X_{n+1}(-2)^{(0,b_{n\!+\!1}^2)}\\
        & \cdots\\
          \cdot\,
        & X_{n-1}(-n\!+\!1)^{(2,-\sum_{k=n}^{n\!+\!1}b_k^{n\!-\!1})}
          X_n(-n\!+\!1)^{(0,b_n^{n\!-\!1})}
          X_{n+1}(-n\!+\!1)^{(0,b_{n\!+\!1}^{n\!-\!1})}\\
          \cdot\,
        & X_n(-n)^{(1,-b_{n\!+\!1}^n)}
          X_{n+1}(-n)^{(0,b_{n\!+\!1}^n)}
\end{aligned}
\end{equation}
where $b_k^i\ge 0$ for all $k,i$. Each element of $\M(\infty)$ may
be written uniquely in this form. Conversely, any element of this
form is an element of $\M(\infty)$.
\end{prop}
\begin{proof}
Given any monomial $M=\prod_{i\in I} \Big({Y_i(-i)}^{(1,a_i^i)}
          \prod_{m=0}^{i-1}{Y_i(-m)}^{(0,a_i^m)}\Big)
          \in\M(\infty)$,
through routine computation using~\eqref{1}, we can obtain the
expression
\begin{equation}~\label{equ:jin19}{\small
\begin{aligned}
\mbox{}\hspace{-4mm}M=& X_1(-1)^{(n,\sum_{i=1}^n a_i^i)}
    X_2(-1)^{(0,-a_1^0)}
    X_3(-1)^{(0,-a_2^0)}
    \cdots
    %X_n(-1)^{(0,-a_{n\!-\!1}^0)}
    X_{n+1}(-1)^{(0,-a_n^0)}\\
    \cdot\,
  & X_2(-2)^{(n\!-\!1,\sum_{i=2}^n a_i^i)}
    X_3(-2)^{(0,-a_1^0\!-\!a_2^1)}
    X_4(-2)^{(0,-a_2^0\!-\!a_3^1)}\cdots
    X_{n+1}(-2)^{(0,-a_{n\!-\!1}^0\!-\!a_n^1)}\\
  & \cdots\\
    \cdot\,
  & X_{n-1}(-n\!+\!1)^{(2,\sum_{i=n\!-\!1}^n a_i^i)}
    X_n(-n\!+\!1)^{(0,-\sum_{i=0}^{n\!-\!2} a_{i+1}^i)}
    X_{n+1}(-n\!+\!1)^{(0,-\sum_{i=0}^{n\!-\!2} a_{i+2}^i)}\\
    \cdot\,
  & X_n(-n)^{(1,a_n^n)}
    X_{n+1}(-n)^{(0,-\sum_{i=0}^{n\!-\!1} a_{i\!+\!1}^i)}.
\end{aligned}}
\end{equation}
Since $M\in\M(\infty)$, from the conditions given
in~(\ref{extend}), we obtain the form given in~\eqref{equ:jin8}.
The element $M_{\infty}=\prod_{i\in I} {Y_i(-i)}^{(1,0)}$ can be
rewritten in the form $\prod_{i\in I}X_i(-i)^{(n-i+1,0)}$.

Conversely, given any monomial of the form \eqref{equ:jin8}, we
have
\begin{equation}~\label{equ:jin49}
\begin{aligned}
\mbox{}\hspace{-4mm}
M=&{Y_1(-1)}^{(1,-\sum_{k=2}^{n\!+\!1}b_k^1+\sum_{k=3}^{n\!+\!1}b_k^2)}
   {Y_1(0)}^{(0,-b_2^1)}\\
   \cdot\, &{Y_2(-2)}^{(1,-\sum_{k=3}^{n\!+\!1}b_k^2
                       +\sum_{k=4}^{n\!+\!1}b_k^3)}
            {Y_2(-1)}^{(0,b_2^1-b_3^2)}
            {Y_2(0)}^{(0,-b_3^1)}\\
   &\ \cdots\\[-2mm]
   \cdot\ &{Y_{n-1}(-n\!+\!1)}^{(1,-\sum_{k=n}^{n\!+\!1}b_k^{n-1}
                                 +\sum_{k=n\!+\!1}^{n\!+\!1}b_k^n)}
           {Y_{n\!-\!1}(-n\!+\!2)}^{(0,b_{n\!-\!1}^{n\!-\!2}-
                                     b_n^{n\!-\!1})}\\
      &\hspace{47mm}\cdots {Y_{n\!-\!1}(-1)}^{(0,b_{n\!-\!1}^1-b_n^2)}
             {Y_{n\!-\!1}(0)}^{(0,-b_n^1)}\\
   \cdot\ &{Y_n(-n)}^{(1,-b_{n\!+\!1}^n)}
           {Y_n(-n\!+\!1)}^{(0,b_n^{n\!-\!1}\!-\!b_{n\!+\!1}^n)}\\
      &\hspace{27mm}\cdots {Y_n(-2)}^{(0,b_n^2\!-\!b_{n\!+\!1}^3)}
             {Y_n(-1)}^{(0,b_n^1\!-\!b_{n\!+\!1}^2)}
             {Y_n(0)}^{(0,-b_{n\!+\!1}^1)}.
\end{aligned}
\end{equation}
It is now straightforward to check that $M\in\M(\infty)$. We have
thus shown that $\M(\infty)$ consists of elements of the
form~\eqref{equ:jin8}.

The uniqueness part may be proved through simple computation.
\end{proof}

\begin{remark}
There are other ways to write each element of $\M(\infty)$ as
products of the terms $X_j(m)^{(u,v)}$. The product
form~\eqref{equ:jin8} was chosen because it allows us to relate
monomials of the set $\M(\infty)$ to tableaux in $\T(\infty)$
directly.
\end{remark}

Now, we translate the Kashiwara operator actions~\eqref{actions},
\eqref{actions2} into a form suitable for the new monomial
expression of $\M(\infty)$.

\begin{lemma}\label{anaysiskashi3}
The set $\M(\infty)$ is closed under the action given below: Fix
element
\begin{equation*}
M=\prod_{i\in
I}\Big(X_i(-i)^{(n-i+1,-\sum_{k=i+1}^{n\!+\!1}b_k^i)}
          \prod_{k=i+1}^{n+1}X_k(-i)^{(0,b_k^i)}\Big)
\in\M(\infty).
\end{equation*}
Consider the following finite ordered sequence
      of some components of $M$:
      \begin{align*}
      &{X_{n+1}(-1)}^{(0,b_{n+1}^{1})},
      {X_{n}(-1)}^{(0,b_{n}^{1})},
      \dots,
      {X_2(-1)}^{(0,b_2^{1})},\\%%%%%%%%%%%%
      &\ {X_{n+1}(-2)}^{(0,b_{n+1}^{2})},
      {X_n(-2)}^{(0,b_n^{2})},
      \dots,
      {X_3(-2)}^{(0,b_3^{2})},\\%%%%%%%%%%%%%
      &\ \ \ \cdots,\\%%%%%%%%%%%%%%%
      &\ \ \ \ {X_{n+1}(-n+1)}^{(0,b_{n+1}^{n-1})},
      {X_n(-n+1)}^{(0,b_n^{n-1})},\\%%%%%%%%%%%%%
      &\ \ \ \ \ {X_{n+1}(-n)}^{(0,b_{n+1}^{n})}.
      \end{align*}%
\begin{enumerate}
\item For $i\in I$, under each component
${X_{i+1}(-m)}^{(0,b_{i+1}^m)}$
      of the above sequence, write $b_{i+1}^m$-many $1$'s, and
      under each ${X_i(-m)}^{(0,b_{i}^m)}$, write $b_i^m$-many $0$'s.
\item From this sequence of $1$'s and $0$'s, successively
      cancel out each $(0,1)$-pair to obtain a sequence of $1$'s
      followed by $0$'s \textup{(}reading from left to right\textup{)}.
      This remaining $1$'s and $0$'s sequence is called
      the \emph{$i$-signature of $M$}.
\item We define
      \begin{equation}\label{111}
      \fit M=M\cdot {X_i(-m)}^{(0,-1)}{X_{i+1}(-m)}^{(0,1)}=
      MA_i(-m)^{-1}
      \end{equation}
      if the component ${X_i(-m)}^{(0,b_{i}^m)}$ corresponds
      to the left-most $0$ of the $i$-signature of $M$ and
      \begin{equation*}
      \eit M=M\cdot{X_i(-m)}^{(0,1)}{X_{i+1}(-m)}^{(0,-1)}=
      MA_i(-m)
      \end{equation*}
      if the component ${X_{i+1}(-m)}^{(0,b_{i+1}^m)}$
      corresponds to the right-most $1$.
\item We define $\eit M=0$ if no $1$ remains and
      \begin{equation}\label{222}
      \fit M=M\cdot {X_i(-i)}^{(0,-1)}{X_{i+1}(-i)}^{(0,1)}=
      MA_i(-i)^{-1}
      \end{equation}
      if no $0$ remains.
\end{enumerate}
\end{lemma}
\begin{proof}
We show that the actions satisfy the following properties:
\begin{equation*}
\fit \M(\infty) \subset \M(\infty), \qquad \eit \M(\infty) \subset
\M(\infty) \cup \{0\} \qquad \text{for all} \ \ i\in I.
\end{equation*}

For $M\in\M(\infty)$, if the $i$-signature of $M$ contains at
least one $0$, then the left-most $0$ of the $i$-signature of $M$
corresponds to a component ${X_i(-m)}^{(0,b_i^{m})}$ for some
$m=1,\dots,i-1$ and the exponent of component
${X_i(-m)}^{(0,b_i^{m})}$ corresponding to the left-most $0$ has
the property $(0,b_i^{m})\ge (0,1)$. Thus the monomial $\fit M=M
{X_i(-m)}^{(0,-1)}{X_{i+1}(-m)}^{(0,1)}$ defined in~\eqref{111} is
contained in $\M(\infty)$.

When the $i$-signature of $M$ contains no $0$, we define $\fit M$ as
in~\eqref{222}. Since the exponents of the components $X_i(-i)$ of
$M$ has the property $\ge (0,1)$, $\fit M$ given in~\eqref{222} also
are in $\M(\infty)$. So the set $\M(\infty)$ is closed under the
above operator $\fit$.

Proof for the statements concerning $\eit$ may be done in a
similar manner.
\end{proof}

\begin{lemma}\label{anaysiskashi2}
The operation given in Lemma~\ref{anaysiskashi3} is just another
expression of the Kashiwara operators given on $\Me$, restricted
to $\M(\infty)$.
\end{lemma}
\begin{proof}
As we can see in equations~\eqref{111} to~\eqref{222}, for each
$M\in\M(\infty)$, $\fit M$ can also be expressed in form
$MA_i(-m)^{-1}$($m=1,2,\dots,i$). To show that this operation is
just another interpretation of the Kashiwara operator $\fit$ given
on $\Me$, restricted to $\M(\infty)$, it is enough to show that
$m_f$ defined in~\eqref{emf} for each $M$ is equal to $-m$ of
$MA_i(-m)^{-1}$ given in equations~\eqref{111} to~\eqref{222}.
Note that we can easily see that from $M$ given by
expression~(\ref{extend}), for each $M\in\M(\infty)$,
$\vphit_i(M)>(0,0)$.

Given a monomial $M\in\M(\infty)$, we can express it in the
following two forms.
\begin{align}
&M=\prod_{i\in
I}\Big(X_i(-i)^{(n-i+1,-\sum_{k=i+1}^{n\!+\!1}b_k^i)}
          \prod_{k=i+1}^{n+1}X_k(-i)^{(0,b_k^i)}\Big)\\
 &\begin{aligned}\label{333}
  =&{Y_1(-1)}^{(1,-\sum_{k=2}^{n\!+\!1}b_k^1+\sum_{k=3}^{n\!+\!1}b_k^2)}
   {Y_1(0)}^{(0,-b_2^1)}\\
   &\ \cdots\\[-2mm]
   \cdot\ &{Y_{n-1}(-n\!+\!1)}^{(1,-\sum_{k=n}^{n\!+\!1}b_k^{n-1}
                                 +\sum_{k=n\!+\!1}^{n\!+\!1}b_k^n)}
           {Y_{n\!-\!1}(-n\!+\!2)}^{(0,b_{n\!-\!1}^{n\!-\!2}-
                                     b_n^{n\!-\!1})}\\
      &\hspace{52mm}\cdots {Y_{n\!-\!1}(-1)}^{(0,b_{n\!-\!1}^1-b_n^2)}
             {Y_{n\!-\!1}(0)}^{(0,-b_n^1)}\\
   \cdot\ &{Y_n(-n)}^{(1,-b_{n\!+\!1}^n)}
           {Y_n(-n\!+\!1)}^{(0,b_n^{n\!-\!1}-b_{n\!+\!1}^n)}
           \cdots {Y_n(-1)}^{(0,b_n^1-b_{n\!+\!1}^2)}
             {Y_n(0)}^{(0,-b_{n\!+\!1}^1)}.
\end{aligned}
\end{align}
If the $i$-signature of $M$ contains at least one $0$ and
${X_i(-m)}^{(0,b_i^{m})}$ is the component corresponding to the
left-most $0$ in the $i$-signature of $M$, then $-m$ is one of
$-1,-2,\dots,-i\!+\!1$ and we can obtain
\begin{equation}
-m=\textup{min}\Bigl\{j\ \Big\vert\
    \textup{max}\big\{ \sum_{k\le j} y_i(k)\mid
    j\in\Z\big\}\Bigr\}=m_f
\end{equation}
where $y_i(k)$ is the exponent of $Y_i(k)$ appearing in $M$ given
by expression~\eqref{333}.

If the $i$-signature of $M$ contains no $0$,
\begin{equation*}
-i=\textup{min}\Big\{j\ \Big\vert\
     \textup{max}\big\{ \sum_{k\le j} y_i(k)\mid j\in\Z\big\}\Big\}=m_f.
\end{equation*}

In all cases, we can confirm that $m_f=-m$, where $-m$ is given
through equations~\eqref{111} to~\eqref{222} stating $\fit M =
MA_i(-m)^{-1}$. Proof for the statements concerning $\eit$ may be
done in a similar manner.
\end{proof}

 From the above lemmas, we obtain the following result.

\begin{prop}\label{zaq}
The set $\M(\infty)$ forms a $\uq(\an)$-subcrystal of $\Me$.
\end{prop}

\begin{figure}
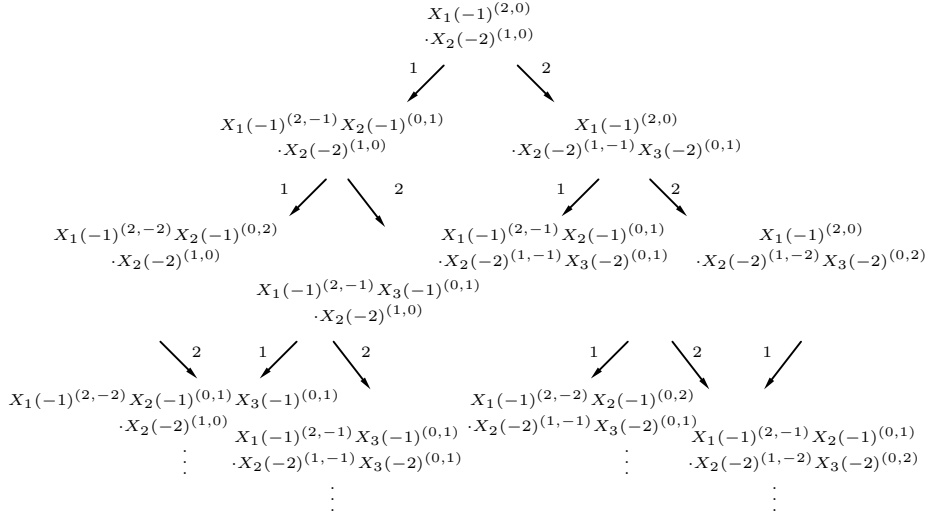

\centering
\begin{texdraw}%
\drawdim in \arrowheadsize l:0.065 w:0.03 \arrowheadtype t:F
\fontsize{6}{6}\selectfont \textref h:C v:C \drawdim em
\setunitscale 1.9
\move(-1 -1.4)\ravec(-1 -1)%%%%%%%%%%%%%%%
\move(1 -1.4)\ravec(1 -1)%%%%%%%%%%%%%%%%%
\move(-3.6 -4.5)\ravec(1 -1.3)%%%%%%%%%%%%
\move(-4.2 -4.5)\ravec(-1 -1)%%%%%%%%%%%%%
\move(3.2 -4.5)\ravec(-1 -1)%%%%%%%%%%%%%%
\move(4.6 -4.5)\ravec(1 -1)%%%%%%%%%%%%%%%
\move(-4 -8.9)\ravec(1 -1.3)%%%%%%%%%%%%%%
\move(-5 -8.9)\ravec(-1 -1)%%%%%%%%%%%%%%%
\move(8.7 -8.9)\ravec(-1 -1.3)%%%%%%%%%%%%
\move(-8.7 -8.9)\ravec(1 -1)%%%%%%%%%%%%%%
\move(4 -8.9)\ravec(-1 -1)%%%%%%%%%%%%%%%%
\move(5.2 -8.9)\ravec(1 -1.3)%%%%%%%%%%%%%
\htext(-8 -12){$\vdots$}%%%%%%%%%%%%%%%%
\htext(4 -12){$\vdots$}%%%%%%%%%%%%%%%%%
\htext(-4 -13){$\vdots$}%%%%%%%%%%%%%%%%
\htext(8 -13){$\vdots$}%%%%%%%%%%%%%%%%%

\htext(0 0){${X_1(-1)}^{(2,0)}$} \htext(0 -0.7){$\cdot
{X_2(-2)}^{(1,0)}$}
%%%%%%%
\htext(-4 -3){${X_1(-1)}^{(2,-1)}{X_2(-1)}^{(0,1)}$} \htext(-4
-3.7){$\cdot {X_2(-2)}^{(1,0)}$}

\htext(4 -3){${X_1(-1)}^{(2,0)}$} \htext(4 -3.7){$\cdot
{X_2(-2)}^{(1,-1)}{X_3(-2)}^{(0,1)}$}
%%%%%%%
\htext(-8.5 -6){${X_1(-1)}^{(2,-2)}{X_2(-1)}^{(0,2)}$} \htext(-8.5
-6.7){$\cdot {X_2(-2)}^{(1,0)}$}

\htext(-3 -7.5){${X_1(-1)}^{(2,-1)}{X_3(-1)}^{(0,1)}$} \htext(-3
-8.2){$\cdot {X_2(-2)}^{(1,0)}$}

\htext(2 -6){${X_1(-1)}^{(2,-1)}{X_2(-1)}^{(0,1)}$} \htext(2
-6.7){$\cdot {X_2(-2)}^{(1,-1)}{X_3(-2)}^{(0,1)}$}

\htext(9 -6){${X_1(-1)}^{(2,0)}$} \htext(9 -6.7){$\cdot
{X_2(-2)}^{(1,-2)}{X_3(-2)}^{(0,2)}$}

\htext(-8.3
-10.5){${X_1(-1)}^{(2,-2)}{X_2(-1)}^{(0,1)}{X_3(-1)}^{(0,1)}$}
\htext(-8.3 -11.2){$\cdot {X_2(-2)}^{(1,0)}$}

\htext(-3.6 -11.5){${X_1(-1)}^{(2,-1)}{X_3(-1)}^{(0,1)}$}
\htext(-3.6 -12.2){$\cdot {X_2(-2)}^{(1,-1)}{X_3(-2)}^{(0,1)}$}

\htext(2.8 -10.5){${X_1(-1)}^{(2,-2)}{X_2(-1)}^{(0,2)}$}
\htext(2.8 -11.2){$\cdot {X_2(-2)}^{(1,-1)}{X_3(-2)}^{(0,1)}$}

\htext(8.8 -11.5){${X_1(-1)}^{(2,-1)}{X_2(-1)}^{(0,1)}$}
\htext(8.8 -12.2){$\cdot {X_2(-2)}^{(1,-2)}{X_3(-2)}^{(0,2)}$}
%%%%%%%%%%%%%%
\move(0 0) \bsegment \htext(-1.8 -1.5){$1$} \htext(1.8 -1.5){$2$}
\esegment \move(0 -2.7) \bsegment \htext(-5.3 -2.1){$1$}
\htext(-2.2 -2.1){$2$} \htext(2.2 -2.1){$1$} \htext(5.3 -2.1){$2$}
\esegment \move(0 -6.2) \bsegment \htext(-7.7 -3){$2$} \htext(-5.9
-3){$1$} \htext(-3.1 -3){$2$} \htext(3.1 -3){$1$} \htext(5.9
-3){$2$} \htext(7.8 -3){$1$} \esegment \move(0 1) %\move(0 -14)
%\drawbb
\end{texdraw}%
\caption{Part of crystal $\M(\infty)$ for type $A_2$}\label{fig:3}
\end{figure}

Figures~\ref{fig:3} illustrates the top part of crystal
$\M(\infty)$ for type $A_2$. It was obtained by applying the
Kashiwara actions introduced in Lemma~\ref{anaysiskashi3} on the
new expression for elements of $\M(\infty)$. Readers may want to
compare this with Figure~\ref{fig:1} and~\ref{fig:2}.

The crystal structure of $\Me$ allows us to obtain more general
results that will be introduced below.

\begin{definition}\label{4.9}
Fix any set of positive integers $p_i$ and any integer $r$. Consider
elements of $\Me$ having the form
\begin{equation}\label{extend2}
\begin{aligned}
      M=& \prod_{i\in I} \Big({Y_i(r-i)}^{(p_i,a_i^i)}
          \prod_{m=0}^{i-1}{Y_i(r-m)}^{(0,a_i^m)}\Big)\\
       =       &{Y_1(r-1)}^{(p_1,a_1^1)}{Y_1(r)}^{(0,a_1^0)}\\
        \cdot\ &{Y_2(r-2)}^{(p_2,a_2^2)}{Y_2(r-1)}^{(0,a_2^1)}
                {Y_2(r)}^{(0,a_2^0)}\\
        & \ \cdots\\
        \cdot\ &{Y_n(r-n)}^{(p_n,a_n^n)}{Y_n(r-n\!+\!1)}^{(0,
                                                    a_n^{n\!-\!1})}
        \cdots {Y_n(r-1)}^{(0,a_n^1)}{Y_n(r)}^{(0,a_n^0)}
\end{aligned}
\end{equation}
and satisfying the same condition given to~\eqref{extend}. In case of
$a_i^j=0$ for all $i,j$, we have
\begin{equation}\label{hi1}
      M=\prod_{i\in I} {Y_i(r-i)}^{(p_i,0)}.
\end{equation}
      We denote by $\M(p_1,\dots,p_n;r;\infty)$
      the set of all monomials of the form~(\ref{extend2})
      and by $M_{(p_1,\dots,p_n;r;\infty)}$ the monomial of~\eqref{hi1}.
\end{definition}

A result similar to Proposition~\ref{equ:jin8} may be obtained for
$\M(p_1,\dots,p_n;r;\infty)$.

\begin{prop}\label{4.10}
Each element of $\M(p_1,\dots,p_n;r;\infty)$ may be written
uniquely in this form
\begin{equation}~\label{equ:jin6}
M=\prod_{i\in I}\Big(X_i(r-i)^{
            (\sum_{k=i}^{n}p_k,-\sum_{k=i+1}^{n\!+\!1}b_k^i)}
  \prod_{k=i+1}^{n+1}X_k(r-i)^{(0,b_k^i)}\Big)
\end{equation}
where $b_k^i\ge 0$ for all $k,i$. Conversely, any element in $\Me$
of this form is an element of $\M(p_1,\dots,p_n;r;\infty)$.
\end{prop}

We believe the readers can easily write down the process for
\emph{change of variable} similar to that given
by~\eqref{equ:jin19} and~\eqref{equ:jin49} for
$\M(p_1,\cdots,p_n;r;\infty)$.

The set $\M(\infty)$ is a special case of this set
$\M(p_1,\dots,p_n;r;\infty)$ which is $r=0$ and $p_i=1$ for all $i\in
I$.

\begin{remark}\label{lem6}
It is possible to obtain results similar to
Lemma~\ref{anaysiskashi3} and~\ref{anaysiskashi2} also for the
case $\M(p_1,\cdots,p_n;r;\infty)$ in a manner similar to that for
$\M(\infty)$. Thus we can state that the set
$\M(p_1,\cdots,p_n;r;\infty)$ forms a subcrystal of $\Me$.
\end{remark}

\begin{prop}\label{cor4}
The set $\M(p_1,\dots,p_n;r;\infty)$ forms a subcrystal of $\Me$
isomorphic to $\M(\infty)$ as $\uq(\an)$-crystal.
\end{prop}
\begin{proof}
As mentioned in Remark~\ref{lem6}, the set
$\M(p_1,\dots,p_n;r;\infty)$ forms a subcrystal of $\Me$. Now let
us show that the crystal $\M(p_1,\dots,p_n;r;\infty)$ is
isomorphic to $\M(\infty)$.

First, we define a canonical map
$\phi:\M(\infty)\rightarrow\M(p_1,\dots,p_n;r;\infty)$ by setting
\begin{equation}
\phi(M)=\prod_{i\in I}\Big(X_i(r-i)^{
            (\sum_{k=i}^{n}p_k,-\sum_{k=i+1}^{n\!+\!1}b_k^i)}
  \prod_{k=i+1}^{n+1}X_k(r-i)^{(0,b_k^i)}\Big)
\end{equation}
for $M$ of the form~\eqref{equ:jin8}. The monomial $M_{\infty}$ of
$\M(\infty)$ is mapped onto the vector
$M_{(p_1,\dots,p_n;r;\infty)}$. It is obvious that this map $\phi$
is well-defined and that it is actually bijective.

Note that the monomial of~\eqref{extend2} is the element of
$\M(p_1,\cdots,p_n;r;\infty)$ corresponding to the
monomial~\eqref{extend} of $\M(\infty)$ under the map $\phi$.

The outputs of the functions $\wt$, $\vphi_i$, and $\veps_i$,
defined in~\eqref{structure1}, \eqref{structure2}, and
\eqref{structure3} do not depend on $r$ or on any fixed positive
integers $p_i$. Using Lemma~\ref{anaysiskashi3} and its
counterpart for $\M(p_1,\dots,p_n;r;\infty)$, we can easily show
that the map $\phi$ commutes with the Kashiwara operators. Hence,
the set $\M(p_1,\dots,p_n;r;\infty)$ forms a subcrystal of $\Me$
isomorphic to $\M(\infty)$.
\end{proof}

\begin{remark}
It should be clear from the proof of Proposition~\ref{cor4}, that
in developing any theory for ${\M(p_1,\cdots,p_n;r;\infty)}$ the
actual values of integer $r$ or $(p_1,\dots,p_n)$ will not be very
important. Arguments made for any set of such values can easily be
adapted to applied to other set of such values. Hence, we shall
concentrate on the theory for $\M(\infty)$ only.
\end{remark}

Now, we will show that $\M(\infty)$ is a new description of
$\B(\infty)$ by giving a crystal isomorphism. Recall that the
crystal $\T(\infty)$ gives a realization of the crystal
$\B(\infty)$(see Theorem~\ref{HL}).

Here is one of our main theorem.

\begin{theorem}\label{mainreal5}
There exists a $\uq(\an)$-crystal isomorphism
\begin{equation*}
\T(\infty) \overset {\sim} \longrightarrow \M(\infty)
\end{equation*}
which maps $T_\infty$ to $M_{\infty}$. It means that the crystal
$\M(\infty)$ is the connected component of $\Me$ containing the
maximal vector $M_{\infty}$ of $\wtt(M_{\infty})=\sum_i(1,0)\La_i$
and is isomorphic to $\B(\infty)$.
\end{theorem}
\begin{proof}
We define a canonical map $\Phi:\T(\infty)\rightarrow\M(\infty)$
by setting, for each tableau $T\in\T(\infty)$ with $i$th ($i\in
I$) row consists of $b^i_j$-many $j$-boxes, for each $i<j\le n+1$,
and some number of $i$-boxes,
\begin{equation}
\Phi(T)=\prod_{i\in
I}\Big(X_i(-i)^{(n-i+1,-\sum_{k=i+1}^{n\!+\!1}b_k^i)}
   \prod_{k=i+1}^{n+1}X_k(-i)^{(0,b_k^i)}\Big)
   \in\M(\infty).
\end{equation}
It is obvious that this map $\Phi$ is well-defined and that it is
actually bijective.

The action of Kashiwara operators on $\M(\infty)$ given in
Lemma~\ref{anaysiskashi3} follows the process for defining it on
$\T(\infty)$. Hence, the map $\Phi$ naturally commutes with the
Kashiwara operators $\fit$ and $\eit$.
\end{proof}

\begin{remark}\label{final2}
 From Proposition~\ref{cor4} and Theorem~\ref{mainreal5}, we can conclude
that the crystal $\M(p_1,\cdots,p_n;r;\infty)$ is also isomorphic to
$\B(\infty)$. It means that for positive integers $p_i$ and integer
$r$, the crystal $\M(p_1,\cdots,p_n;r;\infty)$ is the connected
component of $\Me$ containing the maximal vector
$M_{(p_1,\cdots,p_n;r;\infty)}\in\Me$ of
$\wtt(M_{(p_1,\cdots,p_n;r;\infty)})\\=\sum_{i\in I} (p_i,0)\La_i$
and is isomorphic to $\B(\infty)$.
\end{remark}

\begin{example}
We illustrate the correspondence between $\T(\infty)$ and
$\M(\infty)$ for type $A_3$. Consider a monomial of $\M(\infty)$
\begin{align*}
M=& {Y_1(-1)}^{(1,-5)}{Y_1(0)}^{(0,-3)}\\
  &\cdot {Y_2(-2)}^{(1,-1)}{Y_2(-1)}^{(0,1)}{Y_2(0)}^{(0,0)}\\
  &\cdot {Y_3(-3)}^{(1,-1)}{Y_3(-2)}^{(0,1)}{Y_3(-1)}^{(0,0)}
        {Y_3(0)}^{(0,-4)}.
\end{align*}
It can be expressed as \begin{equation}\label{2}\begin{aligned}
\mbox{}\hspace{-10mm}M=& {X_4(-4)}^{(-1,5)}
    {X_3(-3)}^{(-1,5)}
    {X_2(-2)}^{(-1,5)}
    {X_4(-3)}^{(0,3)}
    {X_3(-2)}^{(0,3)}
    {X_2(-1)}^{(0,3)}\\
  &\cdot {X_4(-4)}^{(-1,1)}
         {X_3(-3)}^{(-1,1)}
         {X_4(-3)}^{(0,-1)}
         {X_3(-2)}^{(0,-1)}\\
  &\cdot {X_4(-4)}^{(-1,1)}
         {X_4(-3)}^{(0,-1)}
         {X_4(-1)}^{(0,4)}
\end{aligned}\end{equation}%
by using the second expression of $Y_i(m)^{(u,v)}$ in~\eqref{1}.
On the other hand, from the equation~\eqref{1}, we also obtain
\begin{equation*}
{X_4(-4)}^{(u,v)}={X_1(-1)}^{(-u,-v)}{X_2(-2)}^{(-u,-v)}
                  {X_3(-3)}^{(-u,-v)}.
\end{equation*}
By applying this equation on~\eqref{2}, we can also be expressed
the monomial $M$ as
\begin{align*}
M=& {X_1(-1)}^{(3,-7)}
    {X_2(-1)}^{(0,3)}
    {X_3(-1)}^{(0,0)}
    {X_4(-1)}^{(0,4)}\\
  &\cdot {X_2(-2)}^{(2,-2)}
         {X_3(-2)}^{(0,2)}
         {X_4(-2)}^{(0,0)}\\
  &\cdot {X_3(-3)}^{(1,-1)}
         {X_4(-3)}^{(0,1)}.
\end{align*}
Actually, we can obtain this expression of $M$
 from~(\ref{equ:jin19}) directly.

Hence we have following marginally large tableau as the image of
$M$ under $\Phi^{-1}$.
\begin{equation*}
\Phi^{-1}(M)= \raisebox{-1.3em}{
\begin{texdraw}%
\drawdim in \arrowheadsize l:0.065 w:0.03 \arrowheadtype t:F
\fontsize{5}{5}\selectfont \textref h:C v:C \drawdim em
\setunitscale 1.6
 \move(-9 3)\rlvec(13 0)
 \move(-9 2)\rlvec(13 0)\rlvec(0 1)
 \move(-9 1)\rlvec(5 0)\rlvec(0 2)
 \move(-9 0)\rlvec(2 0)\rlvec(0 3) \move(0 3)\rlvec(0 -1) \move(-1
3)\rlvec(0 -1) \move(-2 3)\rlvec(0 -1) \move(-3 3)\rlvec(0 -1)
\move(3 3)\rlvec(0 -1) \move(2 3)\rlvec(0 -1) \move(1 3)\rlvec(0
-1) \move(-5 3)\rlvec(0 -2) \move(-6 3)\rlvec(0 -2) \move(-8
3)\rlvec(0 -3) \move(-9 3)\rlvec(0 -3) \htext(-8.5 2.5){$1$}
\htext(-8.5 1.5){$2$} \htext(-8.5 0.5){$3$} \htext(-7.5 2.5){$1$}
\htext(-7.5 1.5){$2$} \htext(-6.5 2.5){$1$} \htext(-6.5 1.5){$2$}
\htext(-5.5 2.5){$1$} \htext(-5.5 1.5){$3$} \htext(-4.5 2.5){$1$}
\htext(-4.5 1.5){$3$} \htext(-3.5 2.5){$1$} \htext(-2.5 2.5){$2$}
\htext(-1.5 2.5){$2$} \htext(-0.5 2.5){$2$} \htext(0.5 2.5){$4$}
\htext(1.5 2.5){$4$} \htext(2.5 2.5){$4$} \htext(3.5 2.5){$4$}
\htext(-7.5 0.5){$4$}
\end{texdraw}}%
\in\T(\infty).
\end{equation*}
\end{example}

%%%%%%%%%%%%%%%%%%%%%%%%%%%%%%%%%%%%%%%%%%%%%%%%%%%%%%%%%%%%%%%%%%%%%%%%%%
\section{Monomial description of crystal $\B(\la)$}%
\label{Extended monomial description of crystals B(la)}

We introduced a realization of crystal $\B(\la)$ given in the
monomial language of Kashiwara and Nakajima through
Theorem~\ref{nakakash} and Corollary~\ref{mla}. Now, we give an
explicit monomial description of $\B(\la)$, for the $\an$-case, i.e.
we will give a concrete listing of elements belonging to the
connected component of $\M_c$(or $\dot\Mec$) containing a certain
maximal vector $M\in\M_c$(or $\dot\Mec$) of weight $\la$. Notice
that there is a natural identification between crystals $\M_c$ and
$\dot\Mec$. We will be working on the Nakajima monomial set $\M_c$
for convenience.

The organization for this section will follow that of the previous
section closely and we take the set $c$ identical to that
of~\eqref{setc} used for the $\B(\infty)$ case. We fix
$\la=l_1\La_1\!+\!\cdots\!+\!l_n\La_n \in {P}^+$.

The set we define below was originally obtained by applying
Kashiwara actions $\fit$($i\in I$) continuously on a maximal vector
$\prod_{i\in I} {Y_i(-i)}^{l_i}\in\M$ of weight $\la$. Actually, we
will be able to confirm later that this choice of starting monomial
allows us to relate monomials of the set defined below to tableaux
in $\B(\la)$ naturally.

\begin{definition}
Consider elements of $\M$ having the form
\begin{equation}\label{eq:jini}
\begin{aligned}
M=& \prod_{i\in I} \Big({Y_i(-i)}^{a_i^i}
    \prod_{m=0}^{i-1}{Y_i(-m)}^{a_i^m}\Big)\\
 =&{Y_1(-1)}^{a_1^1}{Y_1(0)}^{a_1^0}\\
    \cdot\ &{Y_2(-2)}^{a_2^2}{Y_2(-1)}^{a_2^1}
                {Y_2(0)}^{a_2^0}\\
  & \ \cdots\\
    \cdot\ &{Y_n(-n)}^{a_n^n}{Y_n(-n\!+\!1)}^{a_n^{n\!-\!1}}
    \cdots {Y_n(-1)}^{a_n^1}{Y_n(0)}^{a_n^0}
\end{aligned}
\end{equation}
and satisfying the following conditions
\begin{enumerate}
\item $a_i^i\ge 0$, $a_i^0\le 0$, and
      $l_i=\sum_{k=0}^{n-i}a_{i+k}^i-\sum_{k=0}^{i-1}a_{n-i+1+k}^k$
      \ for all $i\in I$,
\item $\sum_{k=0}^ja_{i+k}^k\le 0$ and $\sum_{k=0}^ja_{i+k}^i\ge 0$
      \ for $1\le j\le n\!-\!1$, $1\le i\le n\!-\!j$.
\end{enumerate}
Specifically, in the case of $a_i^j=0$ for all $i>j$, we have
\begin{equation}\label{hi2}
M=\prod_{i\in I} {Y_i(-i)}^{l_i}
 ={Y_1(-1)}^{l_1}{Y_2(-2)}^{l_2}\cdots {Y_n(-n)}^{l_n}.
\end{equation}
We denote by $\M(\la)$ the set of all monomials of the
form~(\ref{eq:jini}) and by $M_{\la}$ the monomial of~\eqref{hi2}.
\end{definition}

In the $A_3$ case, we have
\begin{align*}
M_{\La_1+\La_2}&={Y_1(-1)}{Y_2(-2)},\\
\M(\La_1+\La_2)&=\left\{{\small\begin{aligned}
 &{Y_1(-1)}{Y_2(-2)},\
 {Y_1(0)}^{-1}{Y_2(-2)}{Y_2(-1)},\
 {Y_1(-1)}^{2}{Y_2(-1)}^{-1},\\
 &{Y_2(-2)}{Y_2(0)}^{-1},\
 {Y_1(-1)}{Y_1(0)}^{-1},\\
 &{Y_1(-1)}{Y_2(-1)}^{-1}{Y_2(0)}^{-1},\
 {Y_1(0)}^{-2}{Y_2(-1)},\
 {Y_1(0)}^{-1}{Y_2(0)}^{-1}
\end{aligned}}\right\}.
\end{align*}

\begin{remark}\label{hhh}
Since this set $\M(\la)$ contains a maximal vector $M_{\la}$ of
weight $\la$, by Theorem~\ref{nakakash} it is enough to prove that
the set $\M(\la)$ is closed and connected under Kashiwara
operators~(\ref{f1}) and (\ref{e1}) on $\M$ to show that $\M(\la)$
is a description of the crystal $\B(\la)$.
\end{remark}

For convenience, we shall introduce new expressions for elements
of $\M(\la)$ in terms of the new variable
\begin{equation*}
X_i(m)=Y_i(m)Y_{i-1}(m+1)^{-1}\in\M.
\end{equation*}
These correspond to $X_i(m)^{(0,1)}\in\dot\Me$ introduced in
Definition~\ref{xxx} under the identification between crystals $\M$
and $\dot\Me$ mentioned in Remark~\ref{naka}.

\begin{prop}\label{444}
Consider elements of $\M$ having the form
\begin{equation}~\label{equ:jin3}
\begin{aligned}
M=& \prod_{i\in I}\Big(X_i(-i)^{b_i^i}
    \prod_{k=i+1}^{n+1}X_k(-i)^{b_k^i}\Big)\\
 =& X_1(-1)^{b_1^1}
    X_2(-1)^{b_2^1}
    \cdots
    X_n(-1)^{b_n^1}
    X_{n+1}(-1)^{b_{n\!+\!1}^1}\\
    &\cdot
   X_2(-2)^{b_2^2}
    X_3(-2)^{b_3^2}
    \cdots
    X_{n+1}(-2)^{b_{n\!+\!1}^2}\\
  & \ \cdots\\
    &\cdot
   X_n(-n)^{b_n^n}
    X_{n+1}(-n)^{b_{n\!+\!1}^n}
\end{aligned}
\end{equation}
satisfying the conditions
\begin{enumerate}
\item $b_j^i\ge 0$ \ for all $i,j$, \item
$\sum_{j=i}^{n+1}b_j^i=\sum_{k=i}^{n}l_k$ \ for each $i\in I$,
\item $\sum_{k=0}^jb_{i+k}^i\ge\sum_{k=0}^jb_{i+1+k}^{i+1}$
      \ for $0\le j\le n-1$, $1\le i\le n\!-\!\textup{max}\{1,j\}$.
\end{enumerate}
Each element of $\M(\la)$ may be written uniquely in this form.
Conversely, any element of this form is an element of $\M(\la)$.
\end{prop}
\begin{proof}
Given any monomial $M= \prod_{i\in I} \Big({Y_i(-i)}^{a_i^i}
          \prod_{m=0}^{i-1}{Y_i(-m)}^{a_i^m}\Big)
          \in\M(\la)$,
through routine computation, we can obtain the expression
\begin{equation}~\label{equ:jin2}
\begin{aligned}
\mbox{}\hspace{-4mm}M=& X_1(-1)^{\sum_{i=1}^n a_i^i}
    X_2(-1)^{-a_1^0}
    X_3(-1)^{-a_2^0}
    \cdots
    %X_n(-1)^{(0,-a_{n\!-\!1}^0)}
    X_{n+1}(-1)^{-a_n^0}\\
    \cdot\,
  & X_2(-2)^{\sum_{i=2}^n a_i^i}
    X_3(-2)^{-a_1^0\!-\!a_2^1}
    X_4(-2)^{-a_2^0\!-\!a_3^1}\cdots
    X_{n+1}(-2)^{-a_{n\!-\!1}^0\!-\!a_n^1}\\
  & \ \cdots\\
    \cdot\,
  & X_{n-1}(-n\!+\!1)^{\sum_{i=n\!-\!1}^n a_i^i}
    X_n(-n\!+\!1)^{-\!\sum_{i=0}^{n\!-\!2} a_{i+1}^i}
    X_{n+1}(-n\!+\!1)^{-\!\sum_{i=0}^{n\!-\!2} a_{i+2}^i}\\
    \cdot\,
  & X_n(-n)^{a_n^n}
    X_{n+1}(-n)^{-\sum_{i=0}^{n\!-\!1} a_{i\!+\!1}^i}.
\end{aligned}
\end{equation}
Since $M\in\M(\la)$, from the condition given in~(\ref{eq:jini}),
we obtain the form given in~\eqref{equ:jin3}. The element
$M_{\la}=\prod_{i\in I} {Y_i(-i)}^{l_i}$ can be rewritten in the
form $\prod_{i\in I}X_i(-i)^{\sum_{k=i}^{n}l_k}$.

Conversely, given any monomial of the form \eqref{equ:jin3}, we
have
\begin{equation*}
\begin{aligned}
M=&{Y_1(-1)}^{b_1^1-b_2^2}
   {Y_1(0)}^{-b_2^1}\\
   \cdot\, &{Y_2(-2)}^{b_2^2-b_3^3}
            {Y_2(-1)}^{b_2^1-b_3^2}
            {Y_2(0)}^{-b_3^1}\\
   &\ \cdots\\
   \cdot\ &{Y_{n-1}(-n\!+\!1)}^{b_{n-1}^{n-1}-b_n^n}
           {Y_{n\!-\!1}(-n\!+\!2)}^{b_{n\!-\!1}^{n\!-\!2}-
                                     b_n^{n\!-\!1}}
      \cdots {Y_{n\!-\!1}(-1)}^{b_{n\!-\!1}^1-b_n^2}
             {Y_{n\!-\!1}(0)}^{-b_n^1}\\
   \cdot\ &{Y_n(-n)}^{b_n^n}
           {Y_n(-n\!+\!1)}^{b_n^{n\!-\!1}-b_{n\!+\!1}^n}
      \cdots {Y_n(-1)}^{b_n^1-b_{n\!+\!1}^2}
             {Y_n(0)}^{-b_{n\!+\!1}^1}.
\end{aligned}
\end{equation*}
It is now straightforward to check that $M\in\M(\la)$. We have thus
shown that $\M(\la)$ consists of elements of the
form~\eqref{equ:jin3}.
\end{proof}

Using this new expression for $\M(\la)$, we can easily obtain the
following fact.

\begin{cor}
For $\mu$, $\tau\in {P}^+$, $\M(\mu+\tau)=\M(\mu)\M(\tau)$, where
$\M(\mu)\M(\tau)=\{M\cdot M'\mid M\in \M(\mu),\ M'\in \M(\tau)\}$.
\end{cor}

Now, we translate the Kashiwara operator actions~\eqref{f1},
\eqref{e1} into a form suitable for the new monomial expression of
$\M(\la)$. In a manner similar to the proofs of
Lemmas~\ref{anaysiskashi3} and \ref{anaysiskashi2} for the
$\M(\infty)$ case, we obtain the following result.

\begin{lemma}\label{anaysiskashi6}
The set $\M(\la)$ is closed under the action given below: Fix
element
\begin{equation}
M=\prod_{i\in I}\Big(X_i(-i)^{b_i^i}
  \prod_{k=i+1}^{n+1}X_k(-i)^{b_k^i}\Big)
\in\M(\la).
\end{equation}
Consider the following finite ordered sequence of components of
$M$
      \begin{align*}
      &{X_{n+1}(-1)}^{b_{n+1}^{1}},
      {X_{n}(-1)}^{b_{n}^{1}},
      \dots,
      {X_1(-1)}^{b_1^{1}},\\%%%%%%%%%%%%
      &\ \ {X_{n+1}(-2)}^{b_{n+1}^{2}},
      {X_n(-2)}^{b_n^{2}},
      \dots,
      {X_2(-2)}^{b_2^{2}},\\%%%%%%%%%%%%%
      &\ \ \quad \cdots,\\%%%%%%%%%%%%%%%
      &\ \ \ \ \ \ {X_{n+1}(-n)}^{b_{n+1}^{n}},
      {X_n(-n)}^{b_n^{n}}.
      \end{align*}%
\begin{enumerate}
\item For $i\in I$, under each component
${X_{i+1}(-m)}^{b_{i+1}^m}$
      from the above sequence, write $b_{i+1}^m$-many $1$'s, and
      under each ${X_{i}(-m)}^{b_{i}^m}$, write $b_i^m$-many $0$'s.
\item From this sequence of $1$'s and $0$'s, successively
      cancel out each $(0,1)$-pair to obtain a sequence of $1$'s
      followed by $0$'s \textup{(}reading from left to right\textup{)}.
      This remaining sequence of $1$'s and $0$'s is called
      the \emph{$i$-signature of $M$}.
\item We define
      \begin{equation*}
      \fit M=M{X_i(-m)}^{-1}{X_{i+1}(-m)}=MA_i(-m)^{-1}
      \end{equation*}
      if ${X_{i}(-m)}^{b_{i}^m}$ is the component
      corresponding to the left-most $0$
      appearing in the $i$-signature of $M$.
      And we define
      \begin{equation*}
      \eit M=M{X_i(-m)}{X_{i+1}(-m)}^{-1}=MA_i(-m)
      \end{equation*}
      if ${X_{i+1}(-m)}^{b_{i+1}^m}$ corresponds
      to the right-most $1$.
\item We define $\eit M=0$ if no $1$ remains
      and $\fit M=0$ if no $0$ remains.
\end{enumerate}
\end{lemma}

\begin{lemma}\label{anaysiskashi4}
The operation given in Lemma~\ref{anaysiskashi6} is just another
expression of the Kashiwara operators given on $\M$, restricted to
$\M(\la)$.
\end{lemma}

 From the above lemmas, we obtain the following result immediately.

\begin{prop}\label{closed}
The set $\M(\la)$ forms a $\uq(\an)$-subcrystal of $\M$.
\end{prop}

Figures~\ref{fig:4} and~\ref{fig:5} are the crystal graphs of
$\M(\La_1+\La_2)$ for type $A_2$.

\begin{figure}
\centering
\begin{texdraw}%
\drawdim in \arrowheadsize l:0.065 w:0.03 \arrowheadtype t:F
\fontsize{6}{6}\selectfont \textref h:C v:C \drawdim em
\setunitscale 1.9
\move(2.7 0)\ravec(1.3 0)%%%%%%%%%%%%%%%%
\move(10 0)\ravec(1.3 0)%%%%%%%%%%%%%%%%%
\move(3.2 -2)\ravec(1.3 0)%%%%%%%%%%%%%%%
\move(10 -4)\ravec(1.3 0)%%%%%%%%%%%%%%%%
\move(0.8 -0.7)\ravec(0 -0.7)%%%%%%%%%%%%
\move(7 -2.7)\ravec(0 -0.7)%%%%%%%%%%%%%%
\move(13 -0.7)\ravec(0 -0.7)%%%%%%%%%%%%%
\move(13 -2.7)\ravec(0 -0.7)%%%%%%%%%%%%%
\htext(0.6 0){${Y_1(-1)}{Y_2(-2)}$}
\htext(7 0){${Y_1(0)}^{-1}{Y_2(-2)}{Y_2(-1)}$}
\htext(13.4 0){${Y_2(-2)}{Y_2(0)}^{-1}$}
\htext(0.6 -2){${Y_1(-1)}^{2}{Y_2(-1)}^{-1}$}
\htext(7 -2){${Y_1(-1)}{Y_1(0)}^{-1}$}
\htext(14 -2){${Y_1(-1)}{Y_2(-1)}^{-1}{Y_2(0)}^{-1}$}
\htext(14 -4){${Y_1(0)}^{-1}{Y_2(0)}^{-1}$}
\htext(7 -4){${Y_1(0)}^{-2}{Y_2(-1)}$}
%%%%%%%%%%%%%%
 \htext(3.2 0.3){$1$} \htext(10.6 0.3){$2$}
 \htext(3.6 -1.7){$1$} \htext(10.5 -3.7){$2$}
 \htext(0.4 -0.9){$2$} \htext(13.4 -0.9){$2$}
 \htext(6.7 -2.9){$1$} \htext(13.4 -2.9){$1$}
\end{texdraw}%
\caption{The crystal graph $\M(\La_1+\La_2)$ for type
$A_2$}\label{fig:4}
\end{figure}

\begin{figure}
\centering
\begin{texdraw}%
\drawdim in \arrowheadsize l:0.065 w:0.03 \arrowheadtype t:F
\fontsize{6}{6}\selectfont \textref h:C v:C \drawdim em
\setunitscale 1.9
\move(3.8 0)\ravec(1 0)%%%%%%%%%%%%%%%
\move(9.3 0)\ravec(1 0)%%%%%%%%%%%%%%%
\move(3.8 -2.5)\ravec(1 0)%%%%%%%%%%%%
\move(9.3 -5)\ravec(1 0)%%%%%%%%%%%%%%
\move(1.4 -1.1)\ravec(0 -0.8)%%%%%%%%%
\move(7 -3.7)\ravec(0 -0.8)%%%%%%%%%%%
\move(12.5 -1.1)\ravec(0 -0.8)%%%%%%%%
\move(12.5 -3.7)\ravec(0 -0.8)%%%%%%%%
\htext(1.4 0){${X_1(-1)}^{2}{X_2(-2)}$}
\htext(7 0){${X_1(-1)}{X_2(-1)}$} \htext(7 -0.7){$\cdot
{X_2(-2)}$}
\htext(12.5 0){${X_1(-1)}{X_3(-1)}$} \htext(12.5 -0.7){$\cdot
{X_2(-2)}$}
\htext(1.4 -2.5){${X_1(-1)}^2{X_3(-2)}$}
\htext(7 -2.5){${X_1(-1)}{X_2(-1)}$} \htext(7 -3.2){$\cdot
{X_3(-2)}$}
\htext(12.5 -2.5){${X_1(-1)}{X_3(-1)}$} \htext(12.5 -3.2){$\cdot
{X_3(-2)}$}
\htext(12.5 -5){${X_2(-1)}{X_3(-1)}$} \htext(12.5 -5.7){$\cdot
{X_3(-2)}$}
\htext(7 -5){${X_2(-1)}^2{X_3(-2)}$}
%%%%%%%%%%%%%%
\htext(4.1 0.3){$1$} \htext(9.7 0.3){$2$} \htext(4.1 -2.1){$1$}
\htext(9.7 -4.6){$2$} \htext(1.1 -1.4){$2$} \htext(12.9 -1.4){$2$}
\htext(6.7 -4){$1$} \htext(12.9 -4){$1$}
\end{texdraw}%
\caption{The crystal $\M(\La_1+\La_2)$ for type
$A_2$}\label{fig:5}
\end{figure}

Now, we will show that $\M(\la)$ is a new description of
$\hwc(\la)$. Recall that we identify elements of the highest
weight crystal $\B(\la)$ with semistandard tableaux of
$\la$-shape, for the $A_n$ case. We define a canonical map
$\Psi:\B(\la)\rightarrow\M(\la)$ by setting $\Psi(S)=M$, where $S$
is the semistandard tableau of shape $\la$ with each $i$th row
made up of $b_k^i$-many $k$-blocks for $i\le k\le n+1$ and
\begin{equation*}
M=\prod_{i\in I}\Big(X_i(-i)^{b_i^i}
  \prod_{k=i+1}^{n+1}X_k(-i)^{b_k^i}\Big)
  \in\M(\la).
\end{equation*}
It is obvious that this map $\Psi$ is well-defined and that it is
actually bijective.

The action of Kashiwara operators on the new expressions for
elements of $\M(\la)$ given in Lemma~\ref{anaysiskashi6} follows
the process for defining it on semistandard tableaux of $\B(\la)$.
Hence the map $\Psi$ naturally commutes with the Kashiwara
operators $\fit$ and $\eit$. Other parts of the proof showing that
$\Psi$ is a crystal isomorphism are easy.

\begin{theorem}\label{mainreal}
The bijection $\Psi$ is a $\uq(\an)$-crystal isomorphism.
Therefore, we have a $\uq(\an)$-crystal isomorphism
\begin{equation*}
\hwc(\la) \overset {\sim} \longrightarrow \M(\la).
\end{equation*}
\end{theorem}

\begin{remark}
As we mentioned in Remark~\ref{hhh}, to prove that $\M(\la)$ is a
description of the crystal $\B(\la)$, it is enough to show that
the set $\M(\la)$ is closed and connected under Kashiwara
operators~(\ref{f1}) and (\ref{e1}) on $\M$. In
Proposition~\ref{closed}, the closedness was shown, so now it
remains to show the connectedness. But, in the above theorem, we
proved that $\M(\la)$ is isomorphic to $\B(\la)$ in a direct
manner.

Now we have an explicit monomial description of $\B(\la)$, which
is a connected component in $\M$ of the maximal vector $M_{\la}$
of weight $\la$.
\end{remark}

We can obtain more general results as with the $\M(\infty)$ case.

For any $r\in\Z$, set $\M(r;\la)$ to be the set of all elements in
$\M$ of the form
\begin{equation*}
\prod_{i\in I} \Big({Y_i(r-i)}^{a_i^i}
    \prod_{m=0}^{i-1}{Y_i(r-m)}^{a_i^m}\Big)
\end{equation*}
satisfying the same condition given to~\eqref{eq:jini}.
Specifically, in the case of $a_i^j=0$ for all $i>j$, we set
\begin{equation*}
M_{(r;\la)}=\prod_{i\in I} {Y_i(r-i)}^{l_i}.
\end{equation*}

Each element of the set $\M(r;\la)$ may be written in the form
\begin{equation}\label{prod}
M=\prod_{i\in I}\Big(X_i(r-i)^{b_i^i}
    \prod_{k=i+1}^{n+1}X_k(r-i)^{b_k^i}\Big)
\end{equation}
satisfying the same conditions given to~\eqref{equ:jin3}.
Conversely, any element of this form is an element of $\M(r;\la)$.

The set $\M(\la)$ is a special case of this set $\M(r;\la)$ which
is $r=0$.

Through a process similar to that given in Proposition~\ref{444}
and~\ref{closed} for $\M(\la)$, we can obtain the result that the
set $\M(r;\la)$ becomes a crystal which is isomorphic to
$\M(\la)$. The isomorphism maps $M_{\la}$ of $\M(\la)$ onto the
vector $M_{(r;\la)}$. Thus, the crystal $\M(r;\la)$ is also a
description of $\B(\la)$. Namely, we gave a concrete listing of
elements belonging to the connected component $\M(r;\la)$ of $\M$
containing a maximal vector $M_{(r;\la)}\in\M$ of weight $\la$.

Finally, we can easily obtain the following fact using the
expression~(\ref{prod}) for $\M(r;\la)$.

\begin{cor}
For $\mu$, $\tau\in {P}^+$, $\M(r;\mu+\tau)=\M(r;\mu)\M(r;\tau)$.
\end{cor}

%%%%%%%%%%%%%%%%%%%%%%%%%%%%%%%%%%%%%%%%%%%%%%%%%%%%%%%%%%%%%%%%%%%%%%%%%%
\bibliographystyle{amsplain}

\providecommand{\bysame}{\leavevmode\hbox to3em{\hrulefill}\thinspace}
\providecommand{\MR}{\relax\ifhmode\unskip\space\fi MR }
% \MRhref is called by the amsart/book/proc definition of \MR.
\providecommand{\MRhref}[2]{%
  \href{http://www.ams.org/mathscinet-getitem?mr=#1}{#2}
}
\providecommand{\href}[2]{#2}

%%%%%%%%%%%%%%%%%%%%%%%%%%%%%%%%%%%%%%%%%%%%%%%%%%%%%%%%%%%%%%%%%%%%%%%%%%
\end{document}